\documentclass[12pt,leqno]{article}
\usepackage{amsfonts}
\usepackage{amsmath, amsthm, amsfonts, amssymb, color}
\usepackage{mathrsfs}
\usepackage{url}
\usepackage{color}
\theoremstyle{definition}

\newcommand{\scr}[1]{\mathscr #1}
\definecolor{wco}{rgb}{0.5,0.2,0.3}

\numberwithin{equation}{section} \theoremstyle{remark}

\newcommand{\ua}{\uparrow}

\title{{\bf   Exponential  Ergodicity  for   Non-Dissipative   McKean-Vlasov SDEs }\footnote{Supported in
 part by  the National Key R\&D Program of China (No. 2020YFA0712900) and NNSFC (11771326, 11831014, 11921001).} }
\author{
{\bf    Feng-Yu Wang$^{a),b)}$    }\\
\footnotesize{ Center for Applied Mathematics, Tianjin University, Tianjin 300072, China }  }

\begin{document}
\allowdisplaybreaks
\def\R{\mathbb R}  \def\ff{\frac} \def\ss{\sqrt} \def\B{\mathbf
B}\def\TO{\mathbb T}
\def\I{\mathbb I_{\pp M}}\def\p<{\preceq}
\def\N{\mathbb N} \def\kk{\kappa} \def\m{{\bf m}}
\def\ee{\varepsilon}\def\ddd{D^*}
\def\dd{\delta} \def\DD{\Delta} \def\vv{\varepsilon} \def\rr{\rho}
\def\<{\langle} \def\>{\rangle} \def\GG{\Gamma} \def\gg{\gamma}
  \def\nn{\nabla} \def\pp{\partial} \def\E{\mathbb E}
\def\d{\text{\rm{d}}} \def\bb{\beta} \def\aa{\alpha} \def\D{\scr D}
  \def\si{\sigma} \def\ess{\text{\rm{ess}}}
\def\beg{\begin} \def\beq{\begin{equation}}  \def\F{\scr F}
\def\Ric{{\rm Ric}} \def\Hess{\text{\rm{Hess}}}
\def\e{\text{\rm{e}}} \def\ua{\underline a} \def\OO{\Omega}  \def\oo{\omega}
 \def\tt{\tilde}
\def\cut{\text{\rm{cut}}} \def\P{\mathbb P} \def\ifn{I_n(f^{\bigotimes n})}
\def\C{\scr C}      \def\aaa{\mathbf{r}}     \def\r{r}
\def\gap{\text{\rm{gap}}} \def\prr{\pi_{{\bf m},\varrho}}  \def\r{\mathbf r}
\def\Z{\mathbb Z} \def\vrr{\varrho} \def\ll{\lambda}
\def\L{\scr L}\def\Tt{\tt} \def\TT{\tt}\def\II{\mathbb I}
\def\i{{\rm in}}\def\Sect{{\rm Sect}}  \def\H{\mathbb H}
\def\M{\scr M}\def\Q{\mathbb Q} \def\texto{\text{o}} \def\LL{\Lambda}
\def\Rank{{\rm Rank}} \def\B{\scr B} \def\i{{\rm i}} \def\HR{\hat{\R}^d}
\def\to{\rightarrow}\def\l{\ell}\def\iint{\int}
\def\EE{\scr E}\def\Cut{{\rm Cut}}\def\W{\mathbf W}
\def\A{\scr A} \def\Lip{{\rm Lip}}\def\S{\mathbb S}
\def\BB{\scr B}\def\Ent{{\rm Ent}} \def\i{{\rm i}}\def\itparallel{{\it\parallel}}
\def\g{{\mathbf g}}\def\Sect{{\mathcal Sec}}\def\T{\mathcal T}\def\V{{\bf V}}
\def\PP{{\bf P}}\def\HL{{\bf L}}\def\Id{{\rm Id}}\def\f{{\bf f}}\def\cut{{\rm cut}}
\def\AA{\zeta}

\def\BL{\scr A}

\maketitle

\begin{abstract}  Under Lyapunov and monotone conditions, the exponential ergodicity in the induced Wasserstein quasi-distance is proved for a class of   non-dissipative McKean-Vlasov SDEs, which strengthen  some recent results established under dissipative conditions in long distance. Moreover, when the SDE is order-preserving, the exponential ergodicity is derived in the Wasserstein distance induced by one-dimensional  increasing functions chosen according to the coefficients of    the equation. \end{abstract}

 \vskip 2cm

\section{Introduction}

Consider   the following second order differential operator on $\R^d$:
$$L:=\ff 1 2 \sum_{i,j=1}^d a_{ij} \pp_i\pp_j + \sum_{i=1}^d b_i \pp_i,$$ where
  $a:=(a_{ij})_{1\le i,j\le d}$ is positive definite and    $C^2$-smooth,   $b:=(b_i)_{1\le i\le d}$ is $C^1$-smooth. The Harris theorem says that if  there exists a Lypaunov function $0\le V\in C^2(\R^d)$ with $\lim_{|x|\to\infty} V(x)=\infty$ such that
\beq\label{LY} L V\le C_0-C_1V\end{equation}  holds for some constants $C_0,C_1>0$, then the diffusion process generated by $L$ is exponentially ergodic, see   \cite[Theorem 1.5]{HMS} for a more general assertion, and see  \cite[Theorem 2.1]{EGZ} for an explicit estimate on the exponential convergence rate when $a=I_d,$ the $d\times d$-identity matrix.
A typical example satisfying \eqref{LY}  is  $L=\ff 1 2 \DD + b\cdot\nn$ with $b\in C^1$ such that
 \beq\label{NDDS}  b(x) =-|x|^{p-2} x ,\ \   |x|\ge 1\end{equation}  holds for some $p\ge 1$. It is easy to see that when $p\ge 2,$ this operator is dissipative in long distance,  i.e.
\beq\label{DDS} \<x-y, b(x)-b(y)\>\le C_3 |x-y|-C_4|x-y|^2\end{equation}  holds for some constants $C_3,C_4>0$.
However, when $p\in [1,2),$ it is fully non-dissipative in the sense that
\beq\label{FND} \sup_{|x-y|=r} \<x-y, b(x)-b(y)\>\ge 0,\ \ r\ge 0.\end{equation}
On the other hand, when $p\in (0,1)$,   the diffusion process   is  not exponential ergodic since   the Poincar\'e inequality fails (see for instance   \cite[Corollary 1.4]{W99}).
% and is only sub-exponential ergodic in the sense of \cite{DFG09,RW01}. 
Therefore, in this example,     $p=1$ is critical for the exponential ergodicity.

In this paper, we aim to extend the above mentioned Harris theorem   to McKean-Vlasov SDEs (also called distribution dependent or mean field SDEs), for which the time-marginal distribution $\mu_t$ of the solution satisfies the following nonlinear Fokker-Planck equation on $\scr P$, the space of probability measures on $\R^d$:
\beq\label{NFK} \pp_t \mu_t = L_{\mu_t}^*\mu_t\end{equation} in the sense that $\mu_t$ is continuous in $t$ under the weak topology and
$$\mu_t(f):=\int_{\R^d} f\d\mu_t= \mu_0(f) +\int_0^t \mu_s(L_{\mu_s}f)\d s,\  \ t\ge 0, f\in C_0^\infty(\R^d),$$
where for any $\mu\in \scr P$,   the operator $L_\mu$ is defined by
\beq\label{LM} L_\mu:= \ff 1 2  \sum_{i,j=1}^d a_{ij} \pp_i\pp_j + \sum_{i=1}^d b_i(\cdot,\mu) \pp_i\end{equation}
for the above mentioned $a$ and  a distribution dependent drift
$$b: \R^d\times\scr P\to\R^d.$$ This nonlinear Fokker-Planck equation can be characterized by   the following McKean-Vlasov SDE    on $\R^d$:
 \beq\label{E} \d X_t= b(X_t,\L_{X_t})\d t+ \si(X_t)\d W_t, \end{equation}
 where $\si: [0,\infty)\times \R^d\to \R^d\otimes\R^m$ such that $\si\si^*=a$, $W_t$ is the $m$-dimensional Brownian motion on a complete filtration probability space $(\OO,\{\F_t\}_{t\ge 0},\P)$, and $\L_{\xi}$ is the distribution of a random variable $\xi$.
Indeed,   according to \cite{BR20}, for any solution of \eqref{E1} with $\mu_t:=\L_{X_t}$ satisfying
 \beq\label{ICC}\int_0^T  \mu_t\big(\|a\|+|b(\cdot,\mu_t)|\big)\d t<\infty,\ \ T\ge 0,\end{equation}
 $\mu_t$ solves \eqref{NFK}; while  a solution of \eqref{NFK} satisfying \eqref{ICC}
 coincides with    $\L_{X_t}$ for a weak solution to   \eqref{E}. Therefore, if \eqref{E} is weakly well-posed
 in a subspace $\hat{\scr P}\subset \scr P$ (i.e. for any initial distribution in $\hat{\scr P}$ it has a unique solution with $\L_{X_t}\in  \hat{\scr P}$),  continuous in $t$ under the weak topology,
 the exponential ergodicity of \eqref{NFK} in $\hat {\scr P}$ is equivalent to that of the SDE \eqref{E} with initial distributions in $\hat{\scr P}$.

In recent years,   different types of    exponential ergodicity have been investigated  for solutions to \eqref{NFK} under the dissipative condition \eqref{DDS} in long distance    and that the dependence of $b(x,\mu)$ on $\mu$ is weak enough.
When $a=I_d,$  see \cite{Ma03} for the exponential ergodicity in $\W_2$, 
\cite{V06} for the ergodicity under  the polynomial mixing property for the associated mean-field particle systems, 
\cite{EBB} for the exponential convergence in the total variation norm for Dirac initial measures,   \cite{GLW} for exponential ergodicity in the $``$mean field entropy", \cite{LWZ20} for exponential ergodicity in the $L^1$-Wasserstein distance. See also  \cite{20RWb} for  the exponential ergodicity in  the relative entropy where $a$ may be non-constant. However, as already mentioned above that the condition \eqref{DDS} excludes fully non-dissipative examples  of $b$ satisfying \eqref{NDDS}  for $p\in [1,2)$.
On the other hand,  in this case
the diffusion process generated by $L:= \DD+b\cdot\nn$ is  exponential ergodic according to the Harris theorem, so that in the spirit of stable perturbations, when
$b(x,\mu)=-\nn |x|^p+b_0(x,\mu)$ for $p\in [1,2)$ and large $|x|$ and $b_0$ is small enough, the exponential ergodicity
for \eqref{NFK} with $a=I_d$ should also hold.  This has been confirmed in \cite[Theorem 3.1]{BUT} for $b(x,\mu):=b_0(x)+\vv b_1(x,\mu)$ with small $\vv>0$, where   $b_0$ and $b_1$ satisfies 
$ \<b_0(x),x\>\le - c_1|x|$ for some constant $c_1>0$ and large $|x|$, $ \|b_1\|_\infty\le c_2$ and $$ |b_0(x)-b_0(y)|+|b_1(x,\mu)-b_1(x,\nu)|\le c_2 (|x-y|+\W_2(\mu,\nu))$$
for some constant $c_2>0$ 
and the $L^2$-Wasserstein distance    $\W_2$. 
In this paper, we will prove a general version of such a result for \eqref{NFK}, which includes  non-constant diffusion coefficient $a$ and non $\W_2$-Lipschitz   $b_1(x,\cdot)$, see Example 1.2 below. 

The main idea of the present study is  to decompose $a$  into $a= \ll I_d+\hat\si\hat\si^*$ for some constant $\ll>0$ and Lipschitz continuous $\hat\si$ as in \cite{PW06}, then for the corresponding McKean-Vlasov SDE we adopt the coupling by reflection for the noise with coefficient $\ss\ll I_d$, and the coupling by parallel displacement for the noise with coefficient $\hat\si$.

The coupling by reflection  was applied in \cite{CW94, CW97, CW} to estimate the first eigenvalue on Riemannian manifolds as well as  the spectral gap for elliptic diffusions, and has been developed in the study of SDEs and SPDEs.   Unlike in the study of classical SDEs (or diffusion processes) for which we may let two marginal processes move together after the first meeting time (i.e. coupling time),      in the distribution dependent setting this is no-longer practicable  since  the difference of marginal   distributions  may
separate the marginal processes after  the coupling time. To fix this problem,  after the coupling time we will take the coupling by parallel displacement for all noises, so that the marginal processes will not move too far away each other.

The remainder of the paper is organized as follows. In Section 2, we investigate the exponential ergodicity of \eqref{NFK} under Lyapunov and monotone conditions, which apply to a class of fully non-dissipative models (see Examples 2.1 and 2.2). In Section 3, we prove the exponential ergodicity under the dissipative condition in long time, which extends some existing  results
to non-constant $a$ (see Example 3.1).  Finally,   Section 4 concerns with the exponential ergodicity for order-preserving McKean-Vlasov SDEs.

\section{Under Lyapunov and monotone conditions}

We will consider the following more general version of \eqref{E} where the coefficients may also depend on the time parameter:
\beq\label{E1} \d X_t= b_t(X_t,\L_{X_t})\d t+ \si_t(X_t)\d W_t, \end{equation}
where
$$\si: [0,\infty)\times \R^d\to \R^d\otimes\R^m,\ \ b: [0,\infty)\times\R^d\times\scr P \to \R^d$$
are measurable and $W_t$ is the $m$-dimensional Brownian motion.
Recall that the SDE  \eqref{E1} is called strongly (respectively, weakly) well-posed   in a subspace $\hat{\scr P}\subset\scr P$, if for any $s\ge 0$ and $\F_s$-measurable initial value $X_s$ with $\L_{X_s}\in \hat{\scr P}$
(respectively, any  $\mu\in \hat{\scr P}$),   \eqref{E1} has a unique solution  from time s (respectively, a unique weak solution  with initial distribution $\mu$ from time $s$) such that the time-marginal of the solution is continuous in $\hat{\scr P}$ under the weak topology.  We call \eqref{E1} well-posed if it is both strongly and weakly well-posed.
In this case,  we denote  $P_{s,t}^*\mu=\L_{X_t}$  for  $X_t$ solving \eqref{E1} from time $s$ with  $\L_{X_s}=\mu\in \hat{\scr P}$, so that
$t\mapsto P_{s,t}^*\mu$ is continuous in $t\ge s$ and
\beq\label{SMP} P_{s,t}^*=P_{r,t}^*P_{s,r}^*,\ \ \ 0\le s\le r\le t.  \end{equation}
When $b$ and $a$ do not depend on $t$, we have $P_{s,t}^*= P_{t-s}^*:= P_{0,t-s}^*, t\ge s$.

\subsection{Main result}

For any $t\ge 0$ and $\mu\in \scr P$, consider the second-order differential operator
\beq\label{LM'} L_{t,\mu}:= \ff 1 2  {\rm tr}\{\si_t\si_t^*\nn^2\}+ b_t(\cdot,\mu)\cdot\nn.\end{equation}
For any probability measure $\mu$ and a measurable function $f$, we denote $\mu(f)=\int f \d \mu$ is the integral exists. We assume the following Lyapunov condition. For any positive measurable function $V$ on $\R^d$, let
$$ \scr P_V:=\{\mu\in \scr P: \mu(V)<\infty\}.$$

\beg{enumerate} \item[{$(H_1)$}]  {\bf (Lyapunov) } There exists a function $0\le V\in C^2(V)$ with $\lim_{|x|\to\infty} V(x) =\infty$ and
\beq\label{H11}  \sup_{t\ge 0,\, x\in\R^d} %\bigg\{&\ff{ |\nn V(x)-\nn V(y)|}{|x-y| \{1+V(x)+V(y)\}}\\
%&+
\ff{ |\si_t(x)\nn V(x)|}{1+V(x)} <\infty, \end{equation}
such that  for some    $K_0,K_1\in L^1_{loc}([0,\infty);\R)$
\beq\label{H12} L_{t,\mu} V\le K_0(t)-K_1(t)V,\ \ t\ge 0, \mu\in \scr P_V. \end{equation}
\end{enumerate}
We remark that the existence of invariant probability measure has been studied in \cite{HSS} under an integrated Lyapunov condition weaker than \eqref{H12}, see also \cite{HRW20} for a recent survey on this topic.

For any $l>0$, consider the    class
$$\Psi_l:= \big\{\psi\in C^2([0,l]; [0,\infty)):\ \psi(0)=\psi'(l)=0, \psi'|_{[0,l)}>0 \big\}.$$
For each $\psi\in \Psi_l,$ we extend it to the half line by setting $\psi(r)=\psi(r\land l)$, so that $ \psi'$ is  non-negative and  Lipschitz continuous with compact support, with
\beq\label{HP} c_\psi:= \sup_{r>0} \ff{r\psi'(r)}{\psi(r)} <\infty.\end{equation}
When $\psi''\le 0$, we have $\|\psi'\|_\infty :=\sup |\psi'| = \psi'(0)$ and $c_\psi=\lim_{r\downarrow 0} \ff{r\psi'(r)}{\psi(r)}=1.$

For any  constant $\bb>0$, the  weighted Wasserstein distance (also called transportation cost) is given by
$$\W_{\psi,\bb V} (\mu,\nu):= \inf_{\pi\in \scr C(\mu,\nu)} \int_{\R^d\times\R^d} \psi (|x-y|) \big(1+\bb V(x)+\bb V(y)\big)\pi(\d x,\d y),\ \ \mu,\nu\in \scr P_V.$$
In general, $\W_{\psi,\bb V}$ is only a quasi-distance on $\scr P_V$ as   the triangle inequality may not hold. But it is complete in  the sense that  any $\W_{\psi,\bb V}$-Cauchy sequence in $\scr P_V$   is convergent.
For any $\mu,\nu\in \scr P_V,$ we introduce 
\beq\label{HW}\beg{split}
 \hat W_{\psi,\bb V}(\mu,\nu)&:=\inf_{\pi\in \C(\mu,\nu)} \ff{\int_{\R^d\times\R^d} \psi(|x-y|) (1+\bb V(x)+\bb V(y))  \pi(\d x, \d y)}{\int_{\R^d\times\R^d}  \psi'(|x-y|) (1+\bb V(x)+\bb V(y)) \pi(\d x, \d y)},
 \end{split}\end{equation} which will come naturally from It\^o's formula for the process $$\psi(|X_t-Y_t|) (1+\bb V(X_t)+\bb V(Y_t))$$ for a coupling  $(X_t,Y_t)$ of the SDE. 
 We observe that 
 $$\sup_{\pi\in \C(\mu,\nu)}\int_{\R^d\times\R^d}  \psi'(|x-y|) (1+\bb V(x)+\bb V(y)) \pi(\d x, \d y)\le 1+\bb\mu(V)+\bb \nu(V),$$
  so that  $\hat \W_{\psi,\bb V}\ge \ff{\W_{\psi,\bb V}(\mu,\nu)}{1+\bb\mu(V)+\bb \nu(V)}$.  As shown in Example 1.2 below that in many cases  
  $$  \hat \W_{\psi,\bb V}\ge \ff{c\W_{\psi,\bb V}(\mu,\nu)}{1+\bb[\mu(V) \land   \nu(V)]}$$ holds for some constant $c>0$.  
  
   Moreover,   let $\|\nn f\|_{\infty}$ be the Lipschitz constant of a real function $f$ on $\R^d$. We   need the following  non-degenerate and monotone conditions.

\beg{enumerate} \item[{$(H_2)$}] {\bf (Non-degeneracy)} There exist  $\aa\in L^1_{loc}([0,\infty);(0,\infty))$ and measurable
$$\hat\si: [0,\infty)\times  \R^d\to \R^d\otimes \R^d $$ with $ \int_0^T\|\nn\hat\si_t \|_{\infty}\d t<\infty$ for $T\in (0,\infty)$, such that
\beq\label{A1} a_t(x):=(\si_t\si^*_t)(x)=\aa_t I_d + (\hat\si_t\hat\si^*_t)(x),\ \ t\ge 0, x\in \R^d.\end{equation}
\item[$(H_3)$] {\bf (Monotonicity) }   $b$ is bounded on bounded set in $[0,\infty)\times \R^d\times\scr P_V$. Moreover,    there exist  $l>0$,     $K,\theta,q_l \in L^1_{loc}([0,\infty); [0,\infty))$   and   $\psi\in \Psi_l,$ such that
    \beq\label{A5} 2\aa_t \psi''(r)+ K_t \psi'(r) \le  - q_l(t) \psi(r), \ \ r\in [0,l], t\ge 0,\end{equation}
\beq\label{A33} \beg{split}  & \<b_t(x,\mu)-b_t(y,\nu), x-y\> +\ff 1 2 \|\hat\si_t(x)-\hat\si_t(y)\|_{HS}^2\\
&\le K_t |x-y|^2 +\theta_t |x-y|\hat\W_{\psi,\bb V}(\mu,\nu),\ \ x,y\in \R^d, \mu,\nu \in \scr P_V, t\ge 0.   \end{split}\end{equation}
 \end{enumerate}

\paragraph{Remark 1.1.} (1)  Since $V\ge 0$ with $V(x)\to\infty$ as $|x|\to\infty$,    we have
\beq\label{AA0} \kk_{l,\bb}(t):= \inf_{|x-y|>l} \ff{K_1(t)V(x)+K_1(t) V(y)-2K_0(t)}{\bb^{-1}+ V(x)+V(y)}\in\R,\ \ l>0,\end{equation}
and      $\kk_{l,\bb}(t)>0$ for large enough $l>0$ and $K_1(t)>0$.

(2) Consider the one-dimensional differential operator  $L= 2\ll \ff{\d^2}{\d r^2}+K\ff{\d}{\d r}$ on $[0,l]. $
 In \eqref{A5} one may take $\psi$ to be the first eigenfunction of $L$ with Dirichlet boundary at $0$ and Neumann boundary at $l$. In this case, $q_l>0$ is the first mixed eigenvalue.

(3)   \eqref{H11} and $(H_2)$   imply that
 \beq\label{AA} \beg{split} \aa_{l,\bb}(t):=& c_\psi\sup_{|x-y|\in (0,l)} \bigg\{\aa_t\ff{|\nn V(x)-\nn V(y)| } {|x-y|\{\bb^{-1} + V(x)+V(y)\}}  \\
 & \quad + \ff{|\{\hat\si_t(x)-\hat\si_t(y)\}[(\hat\si_t(\cdot)^*\nn V)(x)+(\hat\si_t(\cdot)^*\nn V)(y)]|}{|x-y|\{\bb^{-1} + V(x)+V(y)\}}\bigg\}<\infty \end{split}\end{equation}
 for any $\bb,l>0$.
 In many cases, we have $\aa_{l,\bb}\downarrow 0$ as $\bb\downarrow 0$. For instance, it is the case when  $V(x)= \e^{|x|^p}$ for $p\in (0,1)$ and large $|x|$,  and   $\hat\si$ is Lipschitz continuous with $\|\hat\si(x)\|\le c(1+ |x|^{q})$
 for some constants $c>0$ and $q\in (0, 1-p),$ or $V(x)=|x|^k$ for some $k>0$ and large $|x|$.

 \

For   $K_0,q_l$, $\kk_{l,\bb}$ and  $\aa_{l,\bb}$ given in $(H_1)$, $(H_3)$, \eqref{AA0} and   \eqref{AA} respectively, let
\beq\label{AA2} \ll_{l,\bb}(t):= \min\big\{\kk_{l,\bb}(t), \ q_l(t)- 2K_0(t)\bb -\aa_{l,\bb}(t)\big\}.\end{equation}
Since $\aa_{l,\bb}(t)\to 0$ as $\bb\to 0,$ and since $\kk_{l,\bb}(t)>0$ for  $K_1(t)>0$ and large $l>0$, when $K_1(t)>0$ we may  take large $l>0$ and small $\bb>0$ such that $\ll_{l,\bb}(t)>0$.
The main result in this part is the following.

 \beg{thm}\label{T1.1} Assume  $(H_1)$-$(H_3),$ with $\psi''\le 0$ when $\hat\si_t(\cdot)$ is non-constant for some $t\ge 0$.    Then the SDE $\eqref{E1}$ is well-posed in $\scr P_V$, and $P_t^*:=P_{0,t}^*$
 satisfies
  \beq\label{EXP1}  \W_{\psi, \bb V}(P_t^*\mu, P_t^*\nu)\le \e^{-\int_0^t\{\ll_{l,\bb} (s) -\theta_s\}\d s}  \W_{\psi, \bb V}(\mu,\nu),\ \ t\ge 0, \mu,\nu\in \scr P_V.\end{equation}
  Consequently, if $(a,b)$ does not depend on $t$  and $\ll_{l,\bb} >\theta,$ then  $P_t^*$ has a unique invariant probability measure
$\bar\mu\in \scr P_V$ such that
\beq\label{EXP2}  \W_{\psi,\bb V}(P_t^*\mu, \bar\mu)\le  \e^{-(\ll_{l,\bb}-\theta ) t}  \W_{\psi, \bb V}(\mu,\bar\mu),\ \ t\ge 0, \mu \in \scr P_V.\end{equation}
  \end{thm}

  \subsection{An example}

  In the following example, the drift $b_0$   is fully non-dissipative in the sense of \eqref{FND}. As mentioned in Introduction that a critical model for the exponential ergodicity is the diffusion process generated by $\DD- (\nn H)\cdot\nn$ with $H(x)= |x|$ for large $|x|$, which is now covered by this example for $p=1$. Moreover, the following example is not covered by \cite{BUT} since \eqref{BUTT} does not hold even for $\hat\si=0$, because $\log\mu(V)$ is   not $\W_2$-Lipschitz continuous in $\mu$.

\paragraph{Example 2.1.}  Let $a=I_d+\hat\si\hat\si^*$ for some Lipschitz continuous matrix valued function $\hat\si$,   $V(x)= \e^{(1+|x|^2)^{p/2}}$ for some $p\in (0,1]$,  and
$$b(x,\mu) := b_0(x) + \vv \Phi(x,\log \mu(V))$$ for some
$\vv\in [0,1),$ $b_0\in C^1(\R^d)$ with $b_0(x)=- |x|^{-p}x$ for $|x|\ge 1$, and $\Phi\in C_b^1(\R^d\times  [0,\infty); \R^d)$. Let
\beq\label{PW} \tt \W_V (\mu,\nu):= \inf_{\pi\in \C(\mu,\nu)} \int_{\R^d\times\R^d} \{1\land |x-y|\} \cdot \{1+V(x)+V(y)\} \pi(\d x,\d y).\end{equation}
Then when $\vv>0$ is small enough, $P_t^*$ has a unique invariant probability measure $\bar\mu\in \scr P_V$,
 and  there exist constants  $c,q>0$ such that
$$ \tt \W_{V}(P_t^*\mu, \bar\mu)\le c \e^{-q t}  \tt \W_{V}(\mu,\bar\mu),\ \ t\ge 0, \mu \in \scr P_V.$$

\beg{proof} It is easy to see  that $(H_1)$ holds for some constants $K_0,K_1>0$,   $(H_2)$ holds for $\aa=1$.
Since $V(x)\to\infty $ as $|x|\to\infty$, we take $l>0$ such that
$$\inf_{|x-y|\ge l} \big\{K_1V(x)+K_1V(y)-2K_0\big\}\ge 1.$$ So,
in \eqref{AA0} the constant $\kk_{l,\bb}>0$ for all $\bb>0$.
Next, take $\psi\in \Psi_l$ such that \eqref{A5} holds for some $q_l>0$, for instance $\psi$ is the first mixed eigenfunction of $2\ff{\d^2}{\d r^2} +K \ff{\d}{\d r}$ on $[0,l]$ with Dirichlet condition at $0$ and Neumann condition at $l$. Then there exists a constant $c_0>0$ such that 
\beq\label{LPS} |V(x)-V(y)|\le c_0 \psi(|x-y|) (V(x)+V(y)),\ \  x,y\in \R^d.\end{equation} 

Next, since for any $\pi\in \C(\mu,\nu)$ we have 
\beg{align*} &\int_{\R^d\times\R^d} \psi'(|x-y|)(1+\bb V(x)+\bb V(y))\pi(\d x,\d y)\\
&\le \|\psi'\|_\infty \int_{\{|x-y|\le l\}} \big\{1+(1+\e)\bb [V(x) \land V(y)]\big\}\pi(\d x,\d y) \\
& \le (2 +\e) [\mu(V)\land \nu(V)],\ \ \bb\in (0,1],\end{align*}
 \eqref{HW} implies 
  $$\hat W_{\psi,\bb V}(\mu,\nu)\ge \ff{\W_{\psi,\bb V}(\mu,\nu)}{(2+\e) [\mu(V)\land\nu(V)]},\ \ \bb\in (0,1].$$
 Combining this with $\Phi\in C_b^1$ and noting that \eqref{LPS} implies 
 $$|\mu(V)-\nu(V)|\le \inf_{\pi\in \C(\mu,\nu)}\int_{\R^d\times\R^d} |V(x)-V(y)|\pi(\d x,\d y)\le c_0\bb^{-1}\W_{\psi,\bb V} (\mu,\nu)$$
 for some constant $c_0>0$, 
 we find a constant $c_1>0$ such that
\beg{align*} &|b(x,\mu)-b(x,\nu)|\le \vv \|\nn \Phi(x,\cdot)\|_\infty |\log \mu(V)-\log \nu(V)|\\
&\le \ff{\vv \|\nn \Phi(x,\cdot)\|_\infty |\mu(V)-\nu(V)|}{\mu(V)\land \nu(V)} 
\le c_1 \vv \bb^{-1} \hat \W_{\psi,\bb V}(\mu,\nu),\ \ \bb\in (0,1].\end{align*} 
Noting that  $\|\nn b_0\|_\infty+\|\nn \Phi\|_\infty+\|\nn\hat\si\|_\infty<\infty$, this implies 
$(H_3)$ holds for some constant $K>0$ and $\theta= c_1 \vv\bb^{-1}, \bb\in (0,1].$

Finally, as observed in Remark 1.1(3) that for the present $V$ we have  $\aa_{l,\bb}\downarrow 0$ as $\bb \downarrow 0. $  Then in \eqref{AA2},  $\ll_{l,\bb}>0$ for small $\bb\in (0,1]$. Therefore,  by Theorem \ref{T1.1},    when $\vv>0$ is small enough, $P_t^*$ has a unique invariant probability measure $\bar\mu\in \scr P_V$, such that
$$\W_{\psi,\bb V} (P_t^*\mu,\bar\mu) \le \e^{-q t} \W_{\psi,\bb V} (\mu,\bar\mu),\ \ t\ge 0$$ holds for some constant $q>0$.
This completes the proof since
$$ C^{-1} \tt \W_V \le \W_{\psi,\bb V} \le C\tt \W_V$$
holds for some constant $C>1$. \end{proof}

\subsection{Proof of Theorem \ref{T1.1}}

Since $\psi(r):=\psi(r\land l)$ for $\psi\in \Psi_l$ is not second order differentiable at $l$, we introduce the following lemma ensuring  It\^o's formula for $\psi$ of a semi-martingale which will be used frequently in the sequence.

\beg{lem}\label{LNN} Let $\xi_t$ be a non-negative  continuous semi-martingale satisfying
$$\d\xi_t\le A_t\d t+ \d M_t$$
for a local martingale $M_t$ and an integrable adapted process $A_t$. Then for any $\psi\in C^1([0,\infty))$ with $\psi'$ non-negative and Lipschitz continuous, we have
$$\d \psi(\xi_t)\le \psi'(\xi_t)A_t \d t + \ff 1 2 \psi''(\xi_t) \d \<M\>_t + \psi'(\xi_t)\d M_t,$$
where $$\psi''(r):= \limsup_{s\downarrow r} \limsup_{\vv\downarrow 0} \ff{\psi'(s+\vv)-\psi'(s)} \vv,\ \ r\ge 0$$ is a bounded measurable function on $[0,\infty)$. \end{lem}

\beg{proof} By restricting before a stopping time, we may and do assume that $\xi_t, \int_0^t A_s\d s$ and $M_t$ are bounded processes.
For any $n\ge 1$, let
$$\psi_n(r)= n \int_0^\infty \psi(r+s) \e^{-ns} \d s,\ \ r\ge 0.$$
Then each $\psi_n$ is $C^\infty$-smooth, with $\psi_n'\ge 0$, $(\psi_n,\psi_n')\to (\psi,\psi')$ locally uniformly,   $\{\|\psi_n''\|_\infty\}_{n\ge 1}$ uniformly bounded, and by   Fatou's lemma,
\beg{align*} &\limsup_{n \to\infty} \psi_n''(r) \le \limsup_{n\to\infty}  \int_0^\infty \limsup_{\vv\downarrow 0} \ff{\psi'(r+s+\vv)-\psi'(r+s)}\vv n\e^{-ns}\d s\\
&\le    \limsup_{s\downarrow 0}   \limsup_{\vv\downarrow 0} \ff{\psi'(r+s+\vv)-\psi'(r+s)}\vv
=\psi''(r),\ \ r\ge 0.\end{align*}
Therefore, by applying It\^o's formula to $\psi_n(\xi_t)$ and letting $n\to \infty$, we finish the proof.
\end{proof}

\paragraph{A. The well-posedness.}

  For any $T>0$ and a subspace $\hat{\scr P}\subset \scr P$, let $C_w([0,T]; \hat{\scr P})$ be the class of all continuous maps from $[0,T]$ to $\hat{\scr P}$ under the weak topology.

\beg{lem}\label{LN1} Assume   that for some $K\in L_{loc}^1([0,\infty); (0,\infty))$
\beq\label{H12*} L_{t,\mu} V(x)\le \AA_t(1+\mu(V)+V(x)),\ \ t\ge 0, x\in\R^d, \mu\in \scr P_V,\end{equation}
\beq\label{H11'} \|\si_t\nn V(x)|\le \AA_t (1+V(x)), \ \ t\ge 0, x\in\R^d,\end{equation}
 \beq\label{H123} \beg{split} &2\<b_t(x,\mu)-b_t(y,\nu),x-y\>^++ \|\si_t(x)-\si_t(y)\|_{HS}^2\\
 &\le \AA_t |x-y|\big\{|x-y|+ \W_{\psi,V}(\mu,\nu)\big\},\ \ t\ge 0, x,y\in \R^d, \mu,\nu\in \scr P_V.\end{split}\end{equation}
 Then   $\eqref{E1} $ is well-posed  for distributions  in $\scr P_V$ with
 \beq\label{EDD} \E V(X_t)\le \e^{2\int_0^T \AA_s\d s} \E V(X_0) \int_0^T \AA_s\e^{2\int_s^T\AA_r\d r}\d s.\end{equation}  \end{lem}

\beg{proof}   It is easy to see that \eqref{EDD} follows from \eqref{H12*} and It\^o's formula. To prove the well-posedness for distributions  in $\scr P_V$, we adopt a fixed point theorem in distributions. For any $T>0$, $\gg:=\L_{X_0}\in \scr P_V,$ and
$$\mu\in \scr P_{T,V}^\gg:=\big\{\mu\in C_w([0,T];\scr P_V):\ \mu_0=\gg\big\},$$     consider the following SDE
\beq\label{E*} \d X_t^\mu= b_t(X_t^\mu,\mu_t)+\si_t(X_t^\mu) \d W_t,\ \ X_0^\mu=X_0, t\in [0,T].\end{equation}
It is well known that the monotone condition  \eqref{A33} in $(H_3)$  implies the well-posedness of this SDE up to life time, while  the Lyapunov condition \eqref{H12*}   implies
$$\sup_{t\in [0,T]} \E [V(X_t^\mu)]<\infty.$$ Then by the continuity of $X_t^\mu$ in $t$ we conclude that
     $$H (\mu)(\cdot):=\L_{X_\cdot^\mu}\in C_w([0,T];\scr P_V).$$
It remains to prove that $H$ has a unique fixed point $\bar\mu\in \scr P_{V,T}$, so that $X_t^{\bar\mu}$ is the unique solution of \eqref{E1} up to time $T$, and by the modified Yamada-Watanabe principle \cite[Lemma 2.1]{HWb}, this also implies the weak well-posedness of \eqref{E1} up to time $T$.

To prove the existence and uniqueness of the fixed point of $H$, we introduce
$$\scr P_{V,T}^{\gg,N}:= \Big\{\mu\in C_w([0,T]; \scr P_V):\ \mu_0=\gg, \sup_{t\in [0,T]} \e^{-N t} \mu_t(V)\le N(1+\gg(V))\Big\},\ \ N\ge 1.$$
Then  as $N\uparrow \infty,$ we have
$\scr P_{V,T}^{\gg, N} \uparrow \scr P_{V,T}^{\gg}$  as $N\uparrow \infty.$ So, it suffices to find $N_0\ge 1$ such that for any $N\ge N_0$,
$H \scr P_{T,V}^{\gg, N}\subset \scr P_{T,V}^{\gg,N}$ and  $H$ has a unique fixed point in $\scr P_{T,V}^{\gg,N}$. We prove this in the following
two steps.

(a) Construction of $N_0.$ Let
$$c:=\e^{\int_0^T \AA_s\d s},\ \
N_0:= 3c.$$ By It\^o's formula and \eqref{H12*}, for any $N\ge N_0$ and $\mu\in \scr P_{T,V}^{\gg,N},$ we have
\beg{align*} &\e^{-Nt}\E V(X_t^\mu)\le \gg(V)\e^{\int_0^t \AA_s\d s-Nt}+ \int_0^t \AA_s\big\{1+N(1+\gg(V))\big\}\e^{\int_s^t \AA_r\d r -N(t-s)}\d s\\
&\le c \gg (V)+ 2c N (1+\gg(V)) \sup_{t\in [0,T]}\int_0^t \e^{-N(t-s)}\d s\le c\gg(V)+ 2c(1+\gg(V))\le N(1+\gg(V)).\end{align*}
So, $H\scr P_{T,V}^{\gg,N}\subset \scr P_{T,V}^{\gg, N}$ for $N\ge N_0.$

(b) Let $N\ge N_0.$ It remains to prove that $H$ is contractive in $\scr P_{T,V}^{\gg,N}$ under
$$\W_{\psi,V,\ll}(\mu,\nu):= \sup_{t\in [0,T]} \e^{-\ll t} \W_{\psi,V}(\mu_t,\nu_t),\ \ \mu,\nu\in\scr P_{T,V}^{\gg,N}$$
for   large $\ll>0$.

For $\mu,\nu\in \scr P_{T,V}^{\gg,N}$, by \eqref{H123}  and the  It\^o-Tanaka formula, we find   $C_0\in L_{loc}^1([0,\infty);(0,\infty))$ such that
$$\d |X_t^\mu-X_t^\nu| \le C_0(t) (\W_{\psi,\bb V}(\mu_t,\nu_t)+|X_t^\mu-X_t^\nu|)\d t +  \Big\<\ff{X_t^\mu-X_t^\nu}{|X_t^\mu-X_t^\nu|}, \big\{\si_t( X_t^\mu)- \si_t(X_t^\nu)\big\}\d W_t\Big\>.$$
Since $\psi\in \Psi_l,$ by extending to the half-line with $\psi(r):= \psi(r\land l)$, we see that   $\psi' $ is non-negative and   Lipschitz continuous.
By   Lemma \ref{LNN}, $\mu,\nu\in \scr P_{V,T}^{\gg,N}$, and noting that $\psi''\le 0$ when $\si_t$ is non-constant for some $t\ge 0$, we find $C_1\in L^1([0,T];(0,\infty))$ such that
\beq\label{XXX} \beg{split} \d \psi(|X_t^\mu-X_t^\nu|)\le & C_1(t)   \psi( |X_t^\mu-X_t^\nu|)+   \W_{\psi,\bb V}(\mu_t,\nu_t)\big\}\d t \\
& +  \psi'( |X_t^\mu-X_t^\nu|)\Big\<\ff{X_t^\mu-X_t^\nu}{|X_t^\mu-X_t^\nu|}, \big\{\si_t(X_t^\mu)- \si_t(X_t^\nu)\big\}\d W_t\Big\>\end{split}\end{equation} holds for $t\in [0,T]$.

On the other hand, by \eqref{H12*} and $\mu,\nu\in \scr P_{V,T}^{\gg,N},$ we find a constant $K(N)>1$ such that
\beg{align*} \d \big\{V(X_t^\mu)+ V(X_t^\nu)\big\} \le &\ \AA_t \big\{1+\mu_t(V)+\nu_t(V) +V(X_t^\mu)+V(X_t^\nu)\big\} \d t\\
 &+ \big\<\si_t(X_t^\mu) \nn V(X_t^\mu)+ \si_t (X_t^\nu) \nn V(X_t^\nu),   \d W_t\big\>\\
\le  K(N) \AA_t \big\{1  +V(X_t^\mu)&+V(X_t^\nu)\big\} \d t
  + \big\<\si_t(X_t^\mu) \nn V(X_t^\mu)+ \si_t(X_t^\nu) \nn V(X_t^\nu),   \d W_t\big\>.\end{align*}
Combining this with \eqref{HP}, \eqref{H11'},  and \eqref{XXX},  we find $C_2\in L^1([0,T];(0,\infty))$ such that
$$\xi_t:= \psi(|X_t^\mu-X_t^\nu|)\big(1+V(X_t^\mu) + V(X_t^\nu)\big)$$ satisfies
$$\d\xi_t  \le C_2(t) \big[  \xi_t  + (1+V(X_t^\mu)+V(X_t^\nu))    \W_{\psi,V}(\mu_t,\nu_t) \big]\d t  +\d M_t,\ \ t\in [0,T]$$ for some  local martingale $M_t$. Since  $H(\mu),H(\nu)\in \scr P_{V,T}^{\gg,N}$ implies
 $$\E V(X_t^\mu)+ \E V(X_t^\nu)\le N(1+\gg(V)) \e^{NT}=:D(N)<\infty,\ \ t\in [0,T],$$ this and $\xi_0=0$  yield
 $$\e^{-\ll t} \E \xi_t\le  \big(1+D(N)\big)\W_{\psi, V,\ll} (\mu,\nu)\int_0^t C_2(s)\e^{ \int_s^t(C_2(r)-\ll)\d r}  \d s,\ \ t\in [0,T],\ll>0.$$
Noting that $\lim_{\ll\to\infty} \sup_{t\in [0,T]} \int_0^t C_2(s)\e^{ \int_s^t(C_2(r)-\ll)\d r}  \d s=0,$
we conclude that when
  $\ll >0$ is large enough,
\beg{align*} \e^{-\ll t} \W_{\psi, V}(H_t(\mu), H_t(\nu))\le   \e^{-\ll t} \E\xi_t
\le \ff 1 2 \W_{\psi, V,\ll}(\mu,\nu),\ \ t\in [0,T].\end{align*}
Therefore,  $H: \scr P_{V,T}^{\gg,N}\to \scr P_{V,T}^{\gg,N}$ is contractive in $\W_{\psi,T,\ll}$ for   large enough $\ll>0$.
 \end{proof}

\paragraph{B. Construction of coupling.}  Simply denote
$$\psi_{\bb V}(x,y):=    \psi(|x-y|) (1+\bb V(x)+\bb V(y)),\ \ x,y\in \R^d.$$ For $s\ge 0$ and $\mu,\nu\in \scr P_V$, let $X_s$ and $Y_s$ be $\F_s$-measurable random variables such that
 \beq\label{WR} \L_{X_s}=P_s^*\mu, \ \ \L_{Y_s}=P_s^* \nu,\ \  \E  \psi_{\bb V}(X_s,Y_s)= \W_{\psi,\bb V}(P_s^*\mu,P_s^*\nu).\end{equation}
 Let $B_t^1$ and $ B_t^2$ be two independent $d$-dimensional Brownian motions and consider the following
 SDE:
 \beq\label{E*A} \d X_t= b_t(X_t, P_t^*\mu)\d t  + \ss{\aa_t}\d B_t^1+ \hat\si_t(X_t)\d B_t^2,\ \ t\ge s.\end{equation}
  By  $(H_3)$, this SDE is well-posed. Indeed, since $b$ is locally bounded, by Girsanov's transform to  the regular SDE
    $$\d X_t=  \ss{\aa_t}\d B_t^1+ \hat\si_t(X_t)\d B_t^2,\ \ t\ge s$$ up to
   the exit time of a large ball, we construct a weak solution to \eqref{E*A} up to the same stopping time. On the other hand, the monotone condition in $(H_3)$ implies the pathwise uniqueness
   of \eqref{E*A}, then the well-posedness is implied by the Yamada-Watanabe principle. Moreover, the Lyapunov condition in $(H_1)$ ensures the non-explosion.
Since by $(H_2)$ and the definition of $P_t^*$ the solution to   the McKean-Vlasov SDE \eqref{E1} is a weak solution to \eqref{E*A}, the weak uniqueness of \eqref{E*A} implies
that  $\L_{X_t}= P_t^*\mu, t\ge s.$

 To construct the   coupling with reflection, let
   $$u(x,y)=   \ff{x-y}{|x-y|},\ \ x\ne y\in \R^d.$$
 We consider the SDE:
 \beq\label{E*1}  \d Y_t= b_t(Y_t, P_t^*\nu)\d t  + \ss{\aa_t}\big\{I_d-2u(X_t,Y_t)\otimes u(X_t,Y_t) 1_{\{t<\tau\}}\big\} \d B_t^1+ \hat\si_t(Y_t)\d B_t^2 \end{equation}
 for $t\ge s$, where
  $$\tau:= \inf\{t\ge s: Y_t=X_t\}$$ is the coupling time.
 Since the coefficients in noises are  Lipschitz continuous in $Y_t \ne X_t$, by the same argument leading to the well-posedness of \eqref{E*A}, we conclude that
 \eqref{E*1} has a unique solution up to the coupling time $\tau$.
 When $t\ge \tau$, the equation of $Y_t$ becomes
 \beq\label{E*1'} \d Y_t= b_t(Y_t, P_t^*\nu)\d t  + \ss{\aa_t} \d B_t^1+ \hat\si_t(Y_t)\d B_t^2,\end{equation} which is well-posed as explained above. Therefore,
 \eqref{E*1} has a unique solution up to life time. On the other hand, the Lyapunov condition in $(H_1)$ implies that the solution is non-explosive, and by the same reason leading to $\L_{X_t}= P_t^*\mu$, we have $\L_{Y_t}= P_t^*\nu.$

 \paragraph{C. Proof of \eqref{EXP1}.}
 By  $(H_3)$ and the It\^o-Tanaka formula for \eqref{E*A} and \eqref{E*1}, we obtain
 \beg{align*} \d |X_t-Y_t| \le &  \big\{ \theta_t  \hat\W_{\psi,\bb V}  (P_t^*\mu,P_t^*\nu)   + K_t |X_t-Y_t| \big\}\d t   \\
 &+ 2\ss{\aa_t}  \<u(X_t,Y_t), \d B_t^1\> + \<u(X_t,Y_t), (\hat\si_t(X_t)-\hat\si_t(Y_t))\d B_t^2\>,\ \ t<\tau.\end{align*}
By Lemma \ref{LNN} and noting that $\psi''\le 0$ when $\hat\si$ is non-constant,      we get
  \beg{align*} &\d \psi(|X_t-Y_t|) \\
  & \le  \big\{\theta_t \psi'(|X_t-Y_t|)  \hat\W_{\psi,\bb V} (P_t^*\mu,P_t^*\nu) + K_t |X_t-Y_t|  \psi'(|X_t-Y_t|)
  + 2\aa_t   \psi''(|X_t-Y_t|)\big\}\d t   \\
& +  \psi'(|X_t-Y_t|)\Big[2 \ss{\aa_t}  \Big\<u(X_t,Y_t), \d B_t^1\Big\>
 +  \Big\<u(X_t,Y_t), (\hat\si_t(X_t)-\hat\si_t(Y_t))\d B_t^2\Big\>\Big],\ \  t<\tau.\end{align*}
Therefore, \eqref{A5} yields
\beq\label{ITP}  \beg{split} \d \psi(|X_t-Y_t|)  \le  & \big\{\theta_t \psi'(|X_t-Y_t|)  \hat\W_{\psi,\bb V} (P_t^*\mu,P_t^*\nu) -q_l(t)   \psi(|X_t-Y_t|)  1_{\{|X_t-Y_t|<l\}}
   \big\}\d t   \\
& +  \psi'(|X_t-Y_t|)\Big[2 \ss{\aa_t}   \Big\<u(X_t,Y_t), \d B_t^1\Big\> \\
&+  \Big\<u(X_t,Y_t), (\hat\si_t(X_t)-\hat\si_t(Y_t))\d B_t^2\Big\>\Big],\ \  t<\tau.\end{split}\end{equation}
By $(H_1)$ and It\^o's formula, we obtain
\beq\label{ITP2} \beg{split} &\d \{V(X_t)+V(Y_t)\}\le \big\{2K_0(t)-K_1(t)V(X_t)-K_1(t)V(Y_t)\}\d t\\
& +\ss{\aa_t} \big\<\nn V(X_t)+ \nn V(Y_t)- 2\<u(X_t,Y_t), \nn V(Y_t)\> u(X_t,Y_t),\ \d B_t^1\big\>\\
&+ \big\<\hat\si_t(X_t)^*\nn V(X_t)+ \hat\si_t(Y_t)^*\nn V(Y_t),\ \d B_t^2\big\>.\end{split}\end{equation}
This together with \eqref{ITP} yields that
$$\phi_t:= \psi_{\bb V}(X_t,Y_t)=\psi(|X_t-Y_t|) \{1+\bb V(X_t)+\bb V(Y_t)\}$$ satisfies
\beq\label{SPP} \beg{split} \d \phi_t \le &\Big\{ \theta_t\psi'(|X_t-Y_t|) \hat \W_{\psi, \bb V}(P_t^*\mu, P_t^*\nu) \big[1+\bb V(X_t)+\bb V(Y_t)\big] -q_l(t) \phi_t 1_{\{|X_t-Y_t|<l\}} \\
&+ \bb \psi(|X_t-Y_t|) \big[2K_0(t)-K_1(t)V(X_t)-K_1(t)V(Y_t)\big]\\
&+ \bb \psi'(|X_t-Y_t|)\Big(\aa_t |\nn V(X_t)-\nn V(Y_t)|\\
& +  \big|\{\hat\si_t(X_t)-\hat\si_t(Y_t)\}[\hat\si_t(X_t)^*\nn V(X_t) + \hat\si_t(Y_t)^*\nn V(Y_t)]\big|\Big)\Big\}\d t \\
&+\d M_t,\ \ t<\tau \end{split}\end{equation}
for some martingale $M_t$.  Combining     \eqref{HP}, \eqref{AA0} and \eqref{AA},  we  derive
\beg{align*} &  \bb \psi(|X_t-Y_t|) \big\{2K_0(t)-K_1(t)V(X_t)-K_1(t)V(Y_t)\big\}\\
&\le 2K_0(t)\bb \phi_t1_{\{|X_t-Y_t|<l\}} - \kk_{l,\bb} (t)\phi_t1_{\{|X_t-Y_t|\ge l\}},\end{align*}  \beg{align*}
&\bb\psi'(|X_t-Y_t|) \Big\{ \aa_t |\nn V(X_t)-\nn V(Y_t)|\\
&\quad + \big|\{\hat\si_t(X_t)-\hat\si_t(Y_t)\}[\hat\si_t(X_t)^*\nn V(X_t) + \hat\si_t(Y_t)^*\nn V(Y_t)]\big|\Big\}\le \aa_{l,\bb}(t) \phi_t1_{\{|X_t-Y_t|<l\}}.\end{align*}
Hence, it follows from \eqref{SPP} that
\beg{align*} \d \phi_t \le &   \theta_t\psi'(|X_t-Y_t|) \hat \W_{\psi,\bb V}(P_t^*\mu, P_t^*\nu) \{1+\bb V(X_t)+\bb V(Y_t)\} \d t\\
& -\big\{[q_l(t)- \aa_{l,\bb}(t)-2K_0(t)\bb] \phi_t 1_{\{|X_t-Y_t|<l\}} + \kk_{l,\bb}(t)\phi_t1_{\{|X_t-Y_t|\ge l\}}\big\}\d t +\d M_t\\
 \le  & \big\{\theta\psi'(|X_t-Y_t|) \hat \W_{\psi,\bb V}(P_t^*\mu, P_t^*\nu) \{1+\bb V(X_t)+\bb V(Y_t)\} -\ll_{l,\bb}(t) \phi_t\big\}\d t+\d M_t,\ \ t<\tau. \end{align*}
 Since $\phi_{t\land \tau}=0$ for $t\ge \tau$, this implies
\beg{align*}& \e^{\int_0^t\ll_{l,\bb}(s)\d s} \E\phi_{t\land\tau} = \E[\phi_{t\land\tau}  \e^{\int_0^{t\land\tau}\ll_{l,\bb} (s)\d s}]
 \le \e^{\int_0^s\ll_{l,\bb} (r)\d r} \E \phi_s \\
 &+ \E\int_s^{t\land\tau} \e^{\int_0^r\ll_{l,\bb}(p)\d p} \theta_t\psi'(|X_t-Y_t|) \hat \W_{\psi,\bb V}(P_r^*\mu, P_r^*\nu)\big\{1+\bb V(X_r)+\bb V(Y_r)\big\}\d r,\ \ t\ge s.\end{align*}
Therefore, for any $t\ge s,$ we have
\beq\label{W1} \beg{split} &\E\phi_{t\land\tau} \le \e^{-\int_s^t\ll_{l,\bb}(r)\d r} \E\phi_s \\
&+  \e^{\int_s^t|\ll_{l,\bb} |(r)\d r} \E \int_s^{t\land\tau}  \theta_r \hat \W_{\psi,\bb V}(P_r^*\mu, P_r^*\nu)\psi'(|X_t-Y_t|) \big\{1+\bb V(X_r)+\bb V(Y_r)\big\}\d r.\end{split}
\end{equation}

On the other hand, for $t\ge \tau$, by It\^o's formula for \eqref{E*A} and \eqref{E*1'}, and applying \eqref{A33},   we find   $C_1\in L_{loc}^1([0,\infty);(0,\infty))$ such that
\beg{align*} \d \psi(|X_t-Y_t|)\le &\big\{C_1(t) \psi(|X_t-Y_t|)  + \theta_t \psi'(|X_t-Y_t|) \hat \W_{\psi,\bb V} (P_t^*\mu, P_t^*\nu) \big\}\d t\\
& + \psi'(|X_t-Y_t|) \<\{\hat\si_t(X_t)-\hat\si_t(Y_t)\}^* u(X_t,Y_t), \d B_t^2\>.\end{align*}
Combining this with \eqref{ITP2}, we find   $C_2\in L_{loc}^1([0,\infty);(0,\infty))$ such that
$$\d\phi_t \le \Big\{C_2(t) \phi_t   + \theta_t  \hat \W_{\psi,\bb V}(P_t^*\mu,P_t^*\nu) \psi'(|X_t-Y_t|) \big\{1+\bb V(X_t)+\bb V(Y_t)\Big\}\d t +\d M_t,\ \ t\ge \tau$$
for some martingale $M_t$. Therefore, for any $t\ge s$, we have $t\land\tau\ge s$ so that 
\beg{align*} &\E\big[1_{\{t>\tau\}}(\phi_t-\phi_{t\land\tau})\big] \\
&\le   \E \int_{t\land\tau}^t  \e^{\int_r^tC_2(p)\d p} \theta_r\hat \W_{\psi,\bb V}(P_r^*\mu, P_r^*\nu) \psi'(|X_r-Y_r|) \big\{1+\bb V(X_r)+\bb V(Y_r)\big\}\d r\\
&\le  \e^{\int_s^tC_2(p)\d p} \E \int_{t\land\tau}^t \theta_r\hat \W_{\psi,\bb V}(P_r^*\mu, P_r^*\nu) \psi'(|X_r-Y_r|) \big\{1+\bb V(X_r)+\bb V(Y_r)\big\}\d r.\end{align*}
This together with \eqref{W1}, \eqref{WR} and \eqref{HW}  yields
\beg{align*} &\E\phi_t = \E\phi_{t\land \tau} + \E\big[1_{\{t>\tau\}}(\phi_t-\phi_{t\land\tau})\big]
 \le \e^{-\int_0^t\ll_{l,\bb}(r)\d r}\E \phi_s \\
 &\qquad +   \e^{\int_s^t(|\ll_{l,\bb}|+ C_2)(r)\d r} \E \int_s^t \theta_r\hat \W_{\psi,\bb V}(P_r^*\mu, P_r^*\nu) \psi'(|X_r-Y_r|)\big\{1+\bb V(X_r)+\bb V(Y_r)\big\}\d r\\
&\le \e^{-\int_s^t\ll_{l,\bb}(r)\d r}\W_{\psi,\bb V} (P_s^*\mu, P_s^*\nu)   +  \e^{\int_s^t(2|\ll_{l,\bb}|+ C_2(r))\d r}  \int_s^t \theta_r\e^{\int_s^r \ll_{l,\bb}(p)\d p} \E \phi_r \d r,\ \ t\ge s,\end{align*}
where the last step follows from the definition of $\hat \W_{\psi,\bb V}$ which implies
$$\hat \W_{\psi,\bb V}(P_r^*\mu, P_r^*\nu) \le \ff{\E\phi_r}{\E[\psi'(|X_r-Y_r|)\{1+\bb V(X_r)+\bb V(Y_r)\}]}.$$
By Gronwall's lemma, we obtain
\beg{align*}\e^{\int_s^t\ll_{l,\bb}(r)\d r} \E \phi_t  \le  \W_{\psi,\bb V} (P_s^*\mu, P_s^*\nu) \exp\bigg[ \e^{\int_s^t   \{2|\ll_{l,\bb}(r)|+ C_2(r)\}\d r} \int_s^t \theta_r\ \d r \bigg],
\ \ t\ge s.\end{align*}
 Thus, for a.e. $s\ge 0,$
 \beg{align*} \ff{\d^+}{\d s} \W_{\psi,\bb V} (P_s^*\mu, P_s^*\nu)&:= \limsup_{t\downarrow s} \ff{\W_{\psi,\bb V} (P_t^*\mu, P_t^*\nu) -\W_{\psi,\bb V} (P_s^*\mu, P_s^*\nu) }{t-s}\\
 &\le  \limsup_{t\downarrow s} \ff{\E \phi_t -\W_{\psi,\bb V} (P_s^*\mu, P_s^*\nu) }{t-s}\\
 &\le -(\ll_{l,\bb}(s) -\theta_s  )\W_{\psi,\bb V} (P_s^*\mu, P_s^*\nu).\end{align*}
 This implies \eqref{EXP1}.

 \paragraph{D. Proof of \eqref{EXP2}.}  Let   $a,b$ be independent of the time parameter and $$\kk:= \ll_{l,\bb} - \theta >0.$$   We intend to show that $P_t^*$ has an invariant probability measure $\bar\mu\in \scr P_V$, so that \eqref{EXP1} implies \eqref{EXP2} and the uniqueness of the invariant probability measure. This can be done as in the proof of \cite[Theorem 3.1(2)]{W18}
 by verifying that $P_t^*\dd_0$ converges in $\scr P_V$ under $\W_{\psi,\bb V}$ as $t\to\infty$, where $\dd_0$ is the Dirac measure at $0\in \R^d$. Precisely, by \eqref{EXP1} and the semigroup property $P_{t+s}^*= P_t^* P_s^*$ for $s,t\ge 0$ due to \eqref{SMP}, we have
 \beg{align*} &\sup_{s\ge 0} \W_{\psi,\bb V} (P_t^*\dd_0, P_{t+s}^*\dd_0) \le \e^{-\kk t} \sup_{s\ge 0} \E^0 \big[\psi(|X_s|) \{1+\bb V(0) +\bb V(X_s)\}\big] \\
 &\le \|\psi\|_\infty \e^{-\kk t} \Big\{1+\bb V(0)+ \sup_{s\ge 0} \E^0 V(X_s)\Big\},\ \ t\ge 0,\end{align*} where $\E^0$ is the expectation taken for the solution to \eqref{E} with $X_0=0$.
 Since \eqref{H12} yields
 $$\sup_{s\ge 0} \E^0 [V(X_s)] \le V(0)+ \ff{K_0}{K_1} <\infty,$$
 we arrive at
 $$\lim_{t\to\infty} \sup_{s\ge 0} \W_{\psi,\bb V} (P_t^*\dd_0, P_{t+s}^*\dd_0)=0,$$
 so that when $t\to\infty,$ $P_t^*\dd_0$ converges to a probability measure $\bar\mu\in\scr P_V$, which is an invariant measure of $P_t^*$.
 Indeed, in this case   the semigroup property and \eqref{EXP1} imply
$$ \W_{\psi,\bb V}(P_s^*\bar\mu,\bar\mu)= \lim_{t\to\infty } \W_{\psi,\bb V}( P_s^*\bar\mu, P_s^*P_t^*\dd_0) \le  \lim_{t\to\infty } \W_{\psi,\bb V}(  \bar\mu,  P_t^*\dd_0) =0,\ \ s\ge 0.$$

\section{Under dissipative condition in long distance}

For any $\psi\in \Psi,$ where
$$\Psi:=\big\{\psi\in C^2([0,\infty)): \psi(0)=0, \psi'>0, r\psi'(r)+r^2(\psi'')^+(r)\le c r\ \text{for\ some\ constant\ }c>0\big\},$$   the  quasi-distance
$$\W_\psi(\mu,\nu):= \inf_{\pi\in \C(\mu,\nu)} \int_{\R^d\times \R^d} \psi(|x-y|)\pi(\d x,\d y)$$
  on the space
$$\scr P_{\psi}:=\big\{\mu\in \scr P: \|\mu\|_\psi:=\mu(\psi(|\cdot|))<\infty\big\}$$
is complete, i.e. a $\W_\psi$-Cauchy sequence in $\scr P_\psi$ converges with resect to $\W_\psi$.  When $\psi$ is concave, $\W_\psi$ satisfies the triangle inequality and is hence a metric on $\scr P_\psi$.

In this part, we do not assume the Lyapunov condition $(H_1)$ but
use the following condition to replace $(H_3)$.

\beg{enumerate}
\item[$(H_3')$] {\bf ($\psi$-Monotonicity) }  Let  $\psi\in \Psi$, $\gg\in C([0,\infty)$ with $\gg(r)\le K r$ for some constant $K>0$ and all $r\ge 0$,   such that
\beq\label{A2E}  2\aa_t \psi''(r) +(\gg\psi')(r)\le -q_t \psi(r),\ \ r\ge 0\end{equation} holds for some $q\in L_{loc}^1([0,\infty); (0,\infty))$.
Moreover, $b$ is locally bounded on $[0,\infty)\times \R^d\times \scr P_\psi$, and   there exists     $\theta \in L_{loc}^1([0,\infty);(0,\infty))$       such that
\beq\label{A3E} \beg{split}  & \<b_t(x,\mu)-b_t(y,\nu), x-y\> +\ff 1 2 \|\hat\si_t(x)-\hat\si_t(y)\|_{HS}^2\\
&\quad \le |x-y|  \big\{\theta_t  \W_\psi (\mu,\nu)  + \gg(|x-y|)\big\},\ \ t\ge 0, x,y\in \R^d, \mu,\nu\in \scr P_\psi.   \end{split}\end{equation}
\end{enumerate}
When $a=I_d$ and
$$b(x,\mu)= b_0(x)+ \int_{\R^d} Z(x,y)\mu(\d y)$$ for a drift $b_0$ and a Lipschitz continuous map $Z:\R^d\times\R^d\to\R^d$,
the exponential convergence of \eqref{E1} is presented  in \cite[Theorems 2.3 and 2.4]{EGZ} under the condition that
$$\<b_0(x)-b_0(y), x-y\>\le \kk(|x-y|) |x-y|^2,\ \ x,y\in \R^d$$
for some function $\kk\in C((0,\infty)) $ with $\int_0^1 r\kk^+(r)\d r<\infty$ and  $\limsup_{r\to\infty } \kk(r)<0$, and that the Lipschitz constant of $Z$ is small enough. It is clear that in this case \eqref{A3E} holds for $\gg(r):=r\kk(r)$ and $\psi(r)$ comparable with $r$, for which we may choose $\psi\in \Psi$ as in \eqref{*OP} below such that \eqref{A2E} holds for $\aa=1$ and some $q>0$. Therefore, this situation is included in Theorem \ref{T3.1} below.

\subsection{Main results and  example}

 \beg{thm}\label{T3.1} Assume  $(H_2)$ and $(H_3')$, with $\psi''\le 0$ if $\hat\si_t(\cdot)$ is non-constant for some $t\ge 0$. Then  $\eqref{E1}$ is well-posed   in $\scr P_\psi$, and $P_t^*$ satisfies
 \beq\label{EXP1'}  \W_\psi(P_t^*\mu, P_t^*\nu)\le  \e^{-\int_0^t\{q_s-\theta_s\|\psi'\|_\infty\}\d s}  \W_\psi(\mu,\nu),\ \ t\ge 0, \mu,\nu \in \scr P_\psi.\end{equation}
 Consequently, $(b_t,a_t)$ does not depend on $t$,
  $q>\theta \|\psi'\|_\infty,$ and
  \beq\label{SUPP} \sup_{t\ge 0} \|P_t^*\dd_0\|_\psi<\infty \end{equation} which is the case when $\psi''\le 0$,    then $P_t^*$ has a unique invariant probability measure
$\bar\mu\in \scr P_\psi$ such that
\beq\label{EXP2'}  \W_\psi(P_t^*\mu, \bar\mu)\le  \e^{-(q-\theta \|\psi'\|_\infty) t}  \W_\psi(\mu,\bar\mu),\ \ t\ge 0, \mu \in \scr P_\psi.\end{equation}
  \end{thm}

  \beg{proof}  By $(H_2)$ and $(H_3')$, the well-posedness follows from the proof of Lemma \ref{LN1} with $\W_\psi$ replacing $\W_{\psi,V}$, and
  the solution satisfies
  \beq\label{*PR} \sup_{t\in [0,T]} \|P_t^*\mu\|_\psi<\infty,\ \ \mu\in \scr P_\psi, T>0.\end{equation}
We omit the details to save space. It remains to prove \eqref{EXP1'} and the existence of the invariant probability measure $\bar \mu$ in the time homogeneous case.

(1) Proof of  \eqref{EXP1'}.  Let $s\ge 0$ and $\mu,\nu\in \scr P_\psi$. We make use the coupling constructed by \eqref{E*A} and \eqref{E*1} for initial values  $(X_s,Y_s)$
 satisfying
 \beq\label{WR'} \L_{X_s}= P_s^*\mu,\ \ \L_{Y_s}= P_s^*\nu,\ \ \W_\psi(P_s^*\mu, P_s^*\nu)= \E\psi(X_s,Y_s).\end{equation}
By the same reason leading to \eqref{ITP},    by $(H_3')$ for $\psi\in \Psi$ with $\psi''\le 0$ when $\hat\si$ is non-constant, we derive
\beq\label{ITP'}  \beg{split} \d \psi(|X_t-Y_t|)  \le & \big\{\theta \psi'(|X_t-Y_t|)   \W_\psi(P_t^*\mu,P_t^*\nu) -q  \psi(|X_t-Y_t|)
   \big\}\d t   \\
& +  \psi'(|X_t-Y_t|)\Big[2 \ss\ll   \Big\<u(X_t,Y_t), \d B_t^1\Big\> \\
&+  \Big\<u(X_t,Y_t), (\hat\si_t(X_t)-\hat\si_t(Y_t))\d B_t^2\Big\>\Big],\ \  t<\tau.\end{split}\end{equation}
By the same argument leading to \eqref{W1}, this implies
\beq\label{W1'} \E \psi(|X_{t\land\tau}-Y_{t\land\tau}|)\le \e^{-q (t-s)}\E \psi(|X_s-Y_s|) +\theta\|\psi'\|_\infty \int_s^{t\land\tau} \W_\psi(P_r^*\mu, P_r^*\nu)\d r,\ \ t\ge s.\end{equation}
On the other hand, when $t\ge \tau$, by $(H_3')$ and applying
It\^o's formula for \eqref{E*A} and \eqref{E*1'},   we find a constant $C>0$ such that
\beg{align*} \d \psi(|X_t-Y_t|)\le & \{C\psi(|X_t-Y_t|) \d t +\theta\|\psi'\|_\infty\W_{\psi} (P_t^*\mu, P_t^*\nu) \big\}\d t\\
& + \psi'(|X_t-Y_t|) \<\{\hat\si_t(X_t)-\hat\si_t(Y_t)\}^* u(X_t,Y_t), \d B_t^2\>.\end{align*}
Thus,
$$\E\big[1_{\{t>\tau\}}\psi(|X_t-Y_t|\big] \le \theta\|\psi'\|_\infty\e^{C(t-s)}  \E \int_{t\land\tau}^t  \W_{\psi}(P_r^*\mu, P_r^*\nu) \d r,
\ \ t\ge s.$$
Combining this with \eqref{W1'} and \eqref{WR'}, we derive
\beg{align*} \W_\psi(P_t^*\mu,P_t^*\nu)&\le \E\psi(|X_t-Y_t|)  = \E\psi(|X_{t\land\tau}- Y_{t\land \tau}|) + \E\big[1_{\{t>\tau\}}\psi(|X_t-Y_t|\big]\\
 &\le \e^{-q(t-s)}\E \psi(|X_s-Y_s|)  + \theta \|\psi'\|_\infty \e^{C(t-s)}  \int_s^t  \W_{\psi}(P_r^*\mu, P_r^*\nu) \d r\\
&=  \e^{-q(t-s)}\W_{\psi} (P_s^*\mu, P_s^*\nu)   + \theta \|\psi'\|_\infty  \e^{C(t-s)}   \int_s^t  \W_{\psi}(P_r^*\mu, P_r^*\nu) \d r,
 \ \ t\ge s.\end{align*}
 Therefore,
 \beg{align*} \ff{\d^+}{\d s} \W_{\psi} (P_s^*\mu, P_s^*\nu)&:= \limsup_{t\downarrow s} \ff{\W_{\psi} (P_t^*\mu, P_t^*\nu) -\W_{\psi} (P_s^*\mu, P_s^*\nu) }{t-s} \\
 &\le -(q -\theta \|\psi'\|_\infty) \W_{\psi} (P_s^*\mu, P_s^*\nu),\ \ s\ge 0.\end{align*}
 This implies \eqref{EXP1'}.

 (2) Existence of $\bar \mu\in \scr P_\psi$.
  Let $(a_t,b_t)$ do not depend on $t$, and $$\ll:=q-\theta\|\psi'\|_\infty>0.$$
  Then \eqref{EXP1'} implies
  $$\W_\psi(P_t^*\dd_0, P_{t+s}^*\dd_0)\le \e^{-\ll t}\W_\psi(\dd_0, P_s^*\dd_0),\ \ t,s\ge 0.$$
  Combining this with \eqref{SPP} we see that as $t\to\infty$, $\{P_t^*\dd_0$ is a $\W_\psi$-Cauchy family whose limit is an invariant probability measure of $P_t^*.$
 It remains to show that \eqref{SUPP} follows from $\psi''\le 0$ and \eqref{EXP1'}. Indeed, in this case $\W_\psi$ satisfies the triangle inequality so that, for $n$ being the integer part of $t>1$,  \eqref{EXP1'} and \eqref{*PR}  imply
\beg{align*} &\|P_t^*\dd_0\|_\psi=\W_\psi(\dd_0, P_t^*\dd_0)\le \sum_{k=0}^{n-1}\W_\psi(P_{k}^*\dd_0, P_{k+1}^*\dd_0) + \W_\psi(P_{n}^*\dd_0, P_t^*\dd_0)\\
 & \le \sum_{k=0}^{n-1} \e^{-\ll k} \|P_1^*\dd_0\|_\psi + \e^{-n\ll}\sup_{s\in [0,1]} \|P_s^*\dd_0\|_\psi\le \Big(\sup_{s\in [0,1]} \|P_s^*\dd_0\|_\psi\Big)\sum_{k=0}^\infty \e^{-\ll k}<\infty.\end{align*}
Therefore, \eqref{SUPP} holds.
\end{proof}

As a consequence of Theorem \ref{T3.1},   we consider the non-dissipative case where $\nn b_t(\cdot,\mu)(x)$ is  positive definite
in a possibly unbounded set but with bounded  $``$one-dimensional puncture mass" in the sense of  \eqref{H23} below.
Let $\scr P_1=\{\mu\in\scr P: \mu(|\cdot|)<\infty\}$ and
 \beg{align*}  S_b(x):= \sup\big\{\<\nn_v b_t(\cdot,\mu)(x), v\>:\ t\ge 0, |v|\le 1, \mu\in \scr P_1\big\},\ \ x\in \R^d.\end{align*}
\beg{enumerate} \item[$(H_3'')$] There exist constants $\theta_0,\theta_1,\theta_2,\aa \ge 0$     such that
\beq\label{H21}  \ff 1 2 \|\hat\si_t(x)-\hat\si_t(y)\|_{HS}^2\le \theta_0|x-y|^2,\ \ t\ge 0, x,y\in\R^d;\end{equation}
\beq\label{H22} S_b(x)\le\theta_1, \ \ |b_t(x,\mu)-b_t(x,\nu)|\le \varphi \W_1(\mu,\nu),\ \ t\ge 0, x\in\R^d, \mu,\nu\in \scr P_1;\end{equation}
\beq\label{H23} \kk:=  \sup_{x,v\in \R^d, |v|=1}\int_\R 1_{\{S_b(x+sv)>-\theta_2\}}\d s <\infty.\end{equation}
\end{enumerate}
Let $\W_1=\W_\psi$ and $\scr P_1= \scr P_\psi$ for $\psi(r)=r$.

\beg{cor} \label{C3.2} Assume $(H_2)$ and $(H_3'')$. Let
\beq\label{G2} \beg{split}& \gg(r):= (\theta_1+\theta_2) \big\{(\kk r^{-1})\land r\big\}-(\theta_2-\theta_0)r,\ \ r\ge 0,\\
&k:= \ff{2\ll}{\int_0^\infty t\, \e^{\ff 1 {2\ll} \int_0^t \gg(u)\d u}\d t }-\ff{\varphi (\theta_2-\theta_0)}{2\ll}\int_0^\infty t \e^{\ff 1 {2\ll}\int_0^t \gg(u)\d u}\d t. \end{split}\end{equation}
Then there exists a constant $c>0$ such that
$$\W_1(P_t^*\mu,P_t^*\nu) \le c\e^{-k t}\, \W_1(\mu,\nu),\ \ t\ge 0, \mu,\nu\in \scr P_1.$$
If $\theta_2>\theta_0$ and
\beq\label{NC} \varphi <  \ff{4\ll^2}{(\theta_2-\theta_2)(\int_0^\infty t\, \e^{\ff 1 {2\ll} \int_0^t \gg(u)\d u}\d t)^2 },\end{equation}
then  $\kk>0$ and $P_t^*$ has a unique invariant probability measure $\bar\mu\in \scr P_1$ satisfying
$$\W_1(P_t^*\mu,\bar\mu) \le c\e^{-k t}\, \W_1(\mu,\bar\mu),\ \ t\ge 0, \mu\in \scr P_1.$$
\end{cor}

\beg{proof} For  $\gg$   in \eqref{G2},  let
\beg{align*} q:= \ff{2\ll}{\int_0^\infty t\, \e^{\ff 1 {2\ll} \int_0^t \gg(u)\d u}\d t },\ \
\theta:= \ff{\varphi (\theta_2-\theta_0)}{2\ll}\int_0^\infty t\e^{\ff 1 {2\ll}\int_0^t \gg(u)\d u}\d t,\end{align*} and take
\beq\label{*OP} \psi(r):= \int_0^r \e^{-\ff 1 {2\ll} \int_0^s \gg(u)\d u} \int_s^\infty t\e^{\ff 1 {2\ll} \int_0^t \gg(u)\d u}\d t,\ \ r\ge 0.\end{equation}
By Theorem \ref{T3.1}, it suffices to verify
\beg{enumerate} \item[(a)] $\psi\in \Psi$ and $\psi''\le 0$;
\item[(b)] there exists a constant $C>1$ such that $C^{-1}\W_\psi\le \W_1\le C\W_\psi;$
\item[(c)]  \eqref{A2E} and \eqref{A3E} hold.
\end{enumerate}

(a) We have $\psi(0)=0, \psi'(r)>0$ and
\beq\label{*1} \psi''(r)= -\ff{\gg(r)}{2\ll} \e^{-\ff 1 {2\ll}\int_0^r \gg(u)\d u}  \int_r^\infty t \e^{\ff 1 {2\ll} \int_0^t\gg(u)\d u}\d t-r,\ \ r\ge 0.\end{equation}
  To prove $\psi\in\Psi$, it suffices to show $\psi''\le 0$. To this end, take
$$r_0:=  \ff{\ss{\kk(\theta_1+\theta_2)}}{\ss{\theta_2-\theta_0}}. $$
It is easy to see that $\gg$ in \eqref{G2} satisfies
\beq\label{*2} \gg|_{[0,r_0]}\ge 0,\ \ \gg|_{(r_0,\infty)}<0.\end{equation}
Combining this with \eqref{*1} we have $\psi''(r)\le 0$ for $r\le r_0.$ On the other hand, for $r>r_0$ we have $\gg(r)<0$ and
$$\ff r{-\gg(r)} = \ff 1 {(\theta_2-\theta_0)r^{1-p}- (\theta_1+\theta_2)  \kk r^{-(1+p)}}$$ is decreasing in $r>r_0$, so that
\beg{align*} &\int_r^\infty t \e^{\ff 1 {2\ll} \int_0^t\gg(u)\d u}\d t =  \int_r^\infty \ff{2\ll t}{\gg(t)}\Big(\ff{\d}{\d t } \e^{\ff 1 {2\ll} \int_0^t\gg(u)\d u}\Big)\d t\\
&=-\ff{2\ll r}{\gg(r)} \e^{\ff 1 {2\ll} \int_0^r\gg(u)\d u} +2\ll  \int_r^\infty \Big(\ff{\d}{\d t} \ff{2\ll t}{-\gg(t)}\Big)   \e^{\ff 1 {2\ll} \int_0^t\gg(u)\d u}\d t\le -\ff{2\ll r}{\gg(r)} \e^{\ff 1 {2\ll} \int_0^r\gg(u)\d u},\ \ r>r_0.\end{align*}
This together with \eqref{*1} yields $\psi''(r)\le 0$ for $r>r_0.$ In conclusion,    $\psi\in \Psi$.

(b) Since $\psi\in \Psi$ with $\psi''\le 0$  implies that $\psi(r)\le \psi'(0)r$ and $\ff{\psi(r)}{r}$ is decreasing in $r>0$, we have
$\W_\psi\le \psi'(0) \W_1$ and
\beq\label{NC2} \beg{split} &\inf_{r>0}\ff{\psi(r)}{r}= \lim_{r\to\infty} \ff {\psi(r)} r = \lim_{r\to\infty}   \psi'(r)  \\
&= \lim_{r\to\infty}\ff{\int_r^\infty t \exp[\ff 1 {2\ll} \int_0^t\gg(u)\d u]\d t}{\exp[\ff 1 {2\ll} \int_0^r\gg(u)\d u]}\\
&= \lim_{r\to\infty}\ff{2\ll r}{-\gg(r)}=\ff {2\ll} {\theta_2-\theta_0}\in (0,\infty).\end{split}\end{equation}
Thus,
$$\ff 1 {\psi'(0)}\W_\psi \le \W_1(\mu,\nu)\le \ff{\theta_2-\theta_0}{2\ll} \W_\psi.$$

(c) By \eqref{G2} we have
$$ 2\ll\psi''(r)+\gg(r)\psi'(r) =-2\ll r,\ \ r\ge 0.$$
Since $\psi(r)\le \psi'(0)r,$ this implies
$$2\ll\psi''(r)+\gg(r)\psi'(r) \le -\ff{2\ll r}{\psi'(0)r} \psi(r)=-q \psi(r),\ \ r\ge 0.$$
Therefore, \eqref{A2E}  holds.

Next, for $x\ne y$, let $v= \ff{x-y}{|x-y|}.$ Then \eqref{H22} implies
\beq\label{G1} \beg{split}& \<b_t(x,\mu)-b_t(y,\nu), x-y\> \\
&= |x-y|\<b_t(x,\mu)-b_t(y,\mu),v\> +|x-y|\<b_t(y,\mu)-b_t(y,\nu), v\>\\
&\le \varphi |x-y|\W_1(\mu,\nu)+ |x-y|\int_0^{|x-y|} S_b(y+s(x-y)) \d s\\
&=\varphi |x-y|\W_1(\mu,\nu)+ \int_0^{|x-y|^2} S_b(y+sv) \d s,\ \ \mu,\nu\in \scr P_1.\end{split}\end{equation}
On the other hand, by \eqref{H22} and \eqref{H23} we obtain
\beg{align*}& \int_0^{|x-y|^2} S_b(y+sv) \d s\le \theta_1 \int_0^{|x-y|^2} 1_{\{S_b(x+sv)>-\theta_2\}}\d s -\theta_2 \int_0^{|x-y|^2} 1_{\{S_b(x+sv)\le -\theta_2\}}\d s\\
&= (\theta_1+\theta_2)  \int_0^{|x-y|^2} 1_{\{S_b(x+sv)>-\theta_2\}}\d s -\theta_2 |x-y|^2\le  (\theta_1+\theta_2)(\kk \land |x-y|^2)-\theta_2|x-y|^2.\end{align*}
Combining this with \eqref{H21} and \eqref{G1}, we derive   \eqref{A3E}.

\end{proof}

To illustrate  Corollary \ref{C3.2}, we consider  the  following nonlinear PDE    for probability density functions $(\rr_t)_{t\ge 0}$ on $\R^d$:
\beq\label{E01} \pp_t \rr_t=    \ff 1 2\Big\{ {\rm div} (a\nn \rr_t) +\sum_{i,j=1}^d \pp_j\big[\rr_t \pp_i a_{ij}\big] \Big\}  +  {\rm div}  \big\{\rr_t \nn (G + W\circledast\rr_t)\big\}, \end{equation}
 where  $G\in C^2(\R^d), W\in C^2(\R^d\times\R^d)$ and
 \beq\label{AOO}  W\circledast\rr_t:=\int_{\R^d}W(\cdot,y) \rr_t(y)\d y.\end{equation}
  According to the correspondence between
 the nonlinear Fokker-Planck equation \eqref{NFK} and the McKean-Vlasov SDE \eqref{E1}, the exponential ergodicity of $\mu_t(\d x):=\rr_t(x)\d x$
  is equivalent to that of $P_t^*$ associated with \eqref{E1} for
 \beq\label{**0} b(x,\mu):= -\nn G(x) - \int_{\R^d} \big\{\nn W(\cdot,z)(x) \big\}   \mu(\d z),\ \ x\in\R^d, \mu\in\scr P_1.\end{equation}
 When  $a= I_d$ and $W$  is symmetric (i.e. $W(x,y)=W(y,x)$),  the exponential ergodicity in $\W_2$ is derived for $\nn^2\ge \ll I_d$ for some $\ll>0$ and a class of $W$ with locally Lipschitz continuous $\nn W$, 
 the exponential ergodicity in the mean field entropy has  been investigated in \cite{GLW} under the dissipative condition in long distance for small enough $\|\nn_x\nn_y W\|_\infty$, see also \cite{CMV} for  a special setting,   \cite{LWZ20} for the exponential ergodicity in $\W_1$,     \cite{EBB} for the exponential convergence in the total variation norm,  and \cite{20RWb} for the exponential ergodicity  in relative entropy.
  In the following example,    we consider the exponential ergodicity  in $\W_1$ for possibly non-constant $a$.
 Indeed, Corollary \ref{C3.2} also applies to granular type equations with non-constant diffusion coefficients.

\paragraph{Example 3.1.}   Let $a$ satisfy $(H_2)$ with $\hat\si$ satisfying $\eqref{H21}$. Consider \eqref{E01} with   $G\in C^2(\R^d)$ and $W\in C^2(\R^d\times\R^d)$ such that
\beq\label{CC1}\beg{split}
& \nn^2_{G+W(\cdot,z)} \ge  \theta_2   1_{\{|\cdot|\ge \ll_0\}} -\theta_1  1_{\{|\cdot|<\ll_0\}},\ \ z\in \R^d,\\
 & \|\nn_x\nn_y W(x,y) \|\le \tt\theta,\ \ x,y\in \R^d\end{split}\end{equation} holds for some constants $ \ll_0,\theta_1,\theta_2>0.$
Then   the assertion in Corollary \ref{C3.2} holds for  $  \kk=4\ll_0$ and  $(P_t^*\mu)(\d x):=\rr_t(x)\d x$, where  $\rr_t$ solves \eqref{E01} with $\rr_0(x)\d x\in \scr P_1. $

\beg{proof}  It is easy to see that  \eqref{CC1} implies \eqref{H22}.
So, it remains to verify that $\kk$ in  \eqref{H23} satisfies  $\kk\le 4 \ll_0$. By the second inequality in \eqref{CC1} we have
$$S_b(x)\le -\theta_21_{\{|x|\ge \ll_0\}} + \theta_11_{\{|x|<\ll_0\}},\ \ x\in \R^d.$$ For  $x,v\in\R^d$ with $|v|=1$, if there exists $s_0\in\R^d$ such that
$|x+s_0v|<\ll_0$, then
$$|x+sv|\ge |s-s_0|-|x+s_0v|> |s-s_0|-\ll_0.$$ so that
$$\{s\in\R: |x+sv|<\ll_0\}\subset (s_0-2\ll_0, s_0+2\ll_0),$$
which implies
$$\kk:=\sup_{x,v\in\R^d, |v|=1} \int_{\{S_b(x+sv)  >-\theta_2\}}\d s \le  \sup_{x,v\in\R^d, |v|=1} \int_{\{|x+sv|<\ll_0\}}\d s \le 4\ll_0.$$
\end{proof}

\section{Order-preserving McKean-Vlasov SDEs}

In this part,  we consider \eqref{E1} with
 \beq\label{O1} \si(x)= {\rm diag} \{\si_1(x_1),\cdots, \si_d(x_d)\},\ \ b_t(x,\mu)= \big(b_1(t,x,\mu),\cdots, b_d(t,x,\mu)\big),\end{equation}  where  $\{\si_i\}_{1\le i\le d} \subset C(\R)$ and
\beq\label{O2}  \beg{split} &b_i(t,x,\mu):=  \bar b_i(x_i) +\int_{\R^d} Z_i(t,x,y)\mu(\d y),\\
& \ \bar b_i\in C(\R), \ \ Z_i\in C([0,\infty)\times\R^d\times\R^d;\R),\ 1\le i\le d.\end{split} \end{equation}
 Then  $X_t=(X_t^1,\cdots, X_t^d)$ for $(X_t^i)_{1\le i\le d}$ solving the SDEs
 \beq\label{E02} \d X_t^i= \big\{\bar b_i(X_t^i) +   \L_{X_t}\big(Z_i(t,X_t,\cdot)\big)\big\}  \d t  +\si_i(X_t^i)\d W_t^i,\ \ 1\le i\le d, \end{equation}
 where  $\mu(f):= \int_{\R^d} f\d\mu$ for a measure $\mu$ on $\R^d$ and a measurable function $f\in L^1(\mu),\ W_t:=(W_t^i)_{1\le i\le d}$ is a $d$-dimensional Brownian motion on a complete filtration probability space $(\OO,\{\F_t\}_{t\ge 0},\P).$

When $Z=0$ (i.e. without interaction), for each $i$, $X_t^i$ is a one-dimensional diffusion process generated by
$$L_i (r)= \bar b_i(r)  \ff{\d }{\d r} + \ff 1 2 \si_i(r)^2 \ff{\d^2}{\d r^2}.$$
Sharp criteria on the exponential ergodicity  have been established for one-dimensional diffusion processes, see for instance \cite{CW}. These  criteria also apply to the diffusion process
generated by $L(x):= \sum_{i=1}^d L_i(x_i)$ as the components are independent one-dimensional diffusion processes.  We will  investigate the exponential ergodicity for the solution to
 \eqref{E02} by making a distribution dependent perturbation to the $L$-diffusion process.

To this end, we take the following class of functions as alternatives  to the first eigenfunction of $L_i, 1\le i\le d:$
\beg{align*} \Phi:= \Big\{\phi=(\phi_1,\cdots,\phi_d): \ &\phi_i\in C^1(\R),  \lim_{|r|\to\infty} |\phi_i(r)|= \infty, \\
&\phi_i'>0 \ \text{is\ locally\ Lipschitz\ continuous}, 1\le i\le d\Big\}.\end{align*}
For any $\phi\in \Phi,$    $\R^d$ is a Polish space under the metric
$$d_\phi(x,y):= |\phi(x)-\phi(y)|_1=\sum_{i=1}^d |\phi_i(x_i)-\phi_i(y_i)|,\ \ x,y\in\R^d,$$  so that
$$\scr P_\phi:= \{\mu\in \scr P: \mu(|\phi|_1)<\infty\}$$ is a Polish space under the Wasserstein distance
$$\W_\phi(\mu,\nu):=\inf_{\pi\in \C(\mu,\nu)} \int_{\R^d\times\R^d} d_\phi(x,y)\pi(\d x,\d y).$$

\subsection{Main result and example}

  \paragraph{(A)} There exists $\phi\in \Phi$  such that the following conditions hold:
\beg{enumerate} \item[$(A_1)$] $\si, V$  and $Z(t,\cdot)$ (uniformly in $t$)  are  locally   Lipschtiz continuous, and there exits a constant $K>0$  such that
\beg{align*}   |\phi_i'\si_i|\le K(1+|\phi_i|),\ \
 L_i \phi_i^2 \le K(1+\phi_i^2). \ \ 1\le i\le d.\end{align*}
\item[$(A_2)$] $\phi_i\in C^2(\R)$   and there exists a constant $q>0$   such that
$$ L_i\phi_i(r) - L_i\phi_i(s)\le -q |\phi_i(r) -\phi_i(s)| ,\ \ -\infty<s\le r<\infty, 1\le i\le d.$$
\item[$(A_3)$] $\sup_{t\ge 0} |Z(t,0,0)|<\infty,$   for each $1\le i\le d$, $Z_i(t,x,y)$ is    increasing in $(x_j)_{j\ne i}$ and $y$, and there exist constants $\theta_1,\theta_2\ge 0$ such that
 \beg{align*} & \sum_{i=1}^d\big|Z_i(t,x,y)\phi_i'(x_i)-Z_i(t,\bar x,\bar y)\phi_i'(\bar x_i)\big| \\
&\le \theta_1 d_\phi(x,\bar x)  +\theta_2 d_\phi(y,\bar y),\ \ t\ge 0, x  \ge \bar x,  y \ge \bar y.\end{align*}
\end{enumerate}
We see that    {\bf (A)} implies the well-posedness of \eqref{E02}   in $\scr P_\phi$, and the solution is order-preserving, i.e.   for any initial values $X_0, Y_0$ with distributions in $\scr P_\phi$ and $\P(X_0\ge Y_0)=1$, we have
 $\P(X_t\ge Y_t,t\ge 0)=1.$

\beg{thm}\label{T2.2} Assume {\bf (A)}.
Then   $\eqref{E02}$ is well-posed and order-preserving for  distributions in $\scr P_\phi$.
Moreover,
\beq\label{XX1} \W_\phi(P_t^*\mu, P_t^*\nu)\le \e^{-(q-\theta_1-\theta_2)t} \W_\phi(\mu,\nu),\ \ t\ge 0, \mu,\nu\in \scr P_\phi.\end{equation}
Consequently, if $q>\theta_1+\theta_2$, then $P_t^*$ has a unique invariant probability measure $\bar\mu\in \scr P_\phi$ such that
\beq\label{XX2} \W_\phi(P_t^*\mu, \bar\mu)\le \e^{-(q-\theta_1-\theta_2)t} \W_\phi(\mu,\bar\mu),\ \ t\ge 0, \mu \in\scr P_\phi.\end{equation}\end{thm}

To illustrate this result, we consider below a simple example which  includes  the one-dimensional diffusion  process generated by
$$L:= \DD-  \nn H$$  with $H(x):=- c|x|$ for some constant $c>0$ and larger $|\cdot|$, since for $\si=\ss 2$ and $\phi(r)={\rm sgn}(r)\e^{\vv |r|} $ for $|r|\ge 1$ we have
$$\bar b(r):=- \ff{q\phi  + \ff{\si^2} 2 \phi'}{\phi'}(r) =-\ff{q+\vv^2}{\vv},\ \ |r|\ge 1.  $$   This is a critical situation for the exponential ergodicity as explained in Introduction.

\paragraph{Example 4.1.}  For each $1\le i\le d$, let $\si_i\in C^1(\R)$   and  $\phi_i\in C^3(\R)$ with $\phi_i'>0$ and $\phi_i(r)={\rm sgn}(r) \e^{\vv |r|}$ for some $\vv>0$ and  $|r|\ge 1$. For  a constant  $q >0$ we take
\beq\label{A10}\bar b_i  = - \ff{q\phi_i + \ff{\si_i^2} 2 \phi_i''}{\phi'_i},\ \ 1\le i\le d.\end{equation}   Moreover,
for  a constant $\aa>0$ and   functions $G_i\in C^1(\R^d\times\R^d)$  increasing in $(x_j)_{j\ne i}$ and $y$ with
\beq\label{A20} \sum_{i=1}^d |G_i(x,y)-G_i(\bar x, \bar y)|\le  d_\phi(x,\bar x)+ d_\phi(y,\bar y),\ \ x,y,\bar x,\bar y\in \R^d,\end{equation}  we take
\beq\label{A30} Z_i(x,y)= \ff{\aa G_i(x,y)}{\phi_i'(x_i)},\ \ x,y\in\R.\end{equation}
Then {\bf (A)} holds for $\theta_1= \theta_2= \aa$.  Consequently, if $\aa<\ff q {2}$ then  $P_t^*$ has a unique invariant probability measure $\bar\mu\in \scr P_\phi$ such that
$$\W_\phi(P_t^*\mu, \bar\mu)\le \e^{-(q-2\aa)t} \W_\phi(\mu,\bar\mu),\ \ t\ge 0, \mu \in \scr P_\phi.$$

\beg{proof} Obviously,   each $Z_i$ is locally Lipschitz continuous with $Z_i(x,y)$  increasing in $(x_j)_{j\ne i}$ and $y$,  and    \eqref{A10} implies
$$L_i\phi_i(r)=\ff{\si_i(r)^2} 2 \phi_i''(r) +\bar b_i(r) \phi_i'(r)= - q \phi_i(r),\ \ r\in \R.$$     Then
$(A_2)$ holds.  Next,  \eqref{A20} and \eqref{A30} yield
\beg{align*}  \sum_{i=1}^d\big|Z_i(x,y)\phi_i'(x_i)-Z_i(\bar x,\bar y)\phi_i'(\bar x_i)\big|\le \sum_{i=1}^d |G_i(x,y)- G_i(\bar x,\bar y)|
 \le   \aa\big\{d_\phi(x, \bar x)+ d_\phi(y,\bar y)\big\},\end{align*}
 so that $(A_3)$ holds  for $\theta_1=\theta_2=\aa.$  Then the desired assertion follows from Theorem \ref{T2.2}.
\end{proof}

 \subsection{Proof of Theorem \ref{T2.2}}

(1) We first prove the well-posedness by using the fixed-point theorem in measures as in the proof of Lemma \ref{LN1}.  Let $X_0$ be $\F_0$-measurable with $\L_{X_0}\in \scr P_\phi$, and let $T>0$. For any $\mu\in C_w([0,T]; \scr P_\phi)$, consider the SDE
$$\d X_t^\mu= b_t(X_t^\mu,\mu_t)+ \si(X_t^\mu)\d W_t,\ \ t\in [0,T], X_0^\mu=X_0,$$ where $b$ and $\si$ are given in \eqref{O1} and \eqref{O2}.
By $(A_1)$, the coefficients of this SDE are locally Lipschitz continuous, so the SDE is well-posed up to the life time $\tau:=\lim_{n\to\infty} \tau_n$, where
$$\tau_n:=\inf\{t\ge 0: |X_t^\mu|\ge n\},\ \ n\ge 1.$$ By $(A_3)$ with $\bar x=\bar y=0$ we obtain
$$\sum_{i=1}^d |Z_i(t,x,y)\phi_i'(x)|\le c_1(1+|\phi(x)|+|\phi(y)|),$$
which together with $(A_1)$   yields
$$\sum_{i=1}^d  L_i\phi_i^2+ 2\sum_{i=1}^d (\phi_i \phi_i')(x_i) \int_{\R^d} Z_i(x,y)  \mu(\d y) \le  c_2\big(1+ |\phi(x)|^2+ |\phi(x)| \mu(|\phi|)\big)  $$
for some constant $c_2>0$.  Then by It\^o's formula, we obtain
$$\d |\phi|^2(X_t^\mu) \le c_2\big(1+ |\phi(X_t^\mu)|^2+  \mu_t(|\phi|)^2 \big)\d t  +2\sum_{i=1}^d (\phi_i\phi_i'\si_i)((X_t^\mu)_i) \d W_t^i.$$
So, letting $\xi_t:=\ss{1+|\phi (X_t^\mu)|^2},$ we derive
$$\d \xi_t \le c_3 \xi_t \d t + \ff 1 {\xi_t} \sum_{i=1}^d (\phi_i\phi_i'\si_i)((X_t^\mu)_i) \d W_t^i$$
for some constant $c_3>0$ depending on $\mu$. Thus,
$$\sup_{t\in [0,T]} (P_t^*\mu)(|\phi|)\le \sup_{t\in [0,T]} \E\xi_t <\infty.$$
This together with the continuity of $X_\cdot^\mu$ yields $H(\mu):=\L_{X_\cdot^\mu}\in C_w([0,T]; \scr P_\phi)$.
So, as explained in the proof of Lemma \ref{LN1}, it remains to show that $H$ is contractive under the metric
$$\W_{\phi,\ll}(\mu,\nu):= \sup_{t\in [0,T] } \e^{-\ll t} \W_\phi(\mu_t,\nu_t),\ \ \mu,\nu\in C_w([0,T];\scr P_\phi)$$ for large $\ll>0$.

For  $\mu^1,\mu^2\in C_w([0,T];\scr P_\phi)$, we choose random variables $\eta^1_\cdot, \eta^2_\cdot$ on $C([0,T];\R^d)$ such that $\L_{\eta_\cdot^i}=\mu^i, i=1,2.$
Let
$$\bar\mu_t:=\L_{\eta^1_t\lor\eta^2_t}, \ \ \hat \mu_t:= \L_{\eta^1_t\land\eta^2_t},\ \ t\in [0,T].$$
Then $\bar\mu_t\ge \mu_t^i\ge \hat\mu_t  $ in the sense
$$\bar\mu_t(f)\ge\mu_t^i(f)\ge  \hat\mu_t(f),\ \ t\in [0,T], f\in \scr M_b(\R^d),  i=1,2$$  where  $\scr M_b(\R^d)$ is the class of all bounded increasing functions on $\R^d$.
Combining this with $(A_3)$, we conclude that
\beq\label{CMP0} b_i(t,x,\hat \mu_t) \le b_i(t,y,\mu_t^1), b_i(t,y,\mu_t^2) \le b_i(t,z,\bar\mu_t),\ \ t\ge 0\end{equation}
holds for $1\le i\le d$ and $x,y,z\in \R^d$ with $ x_i=y_i=z_i$ and $ x_j\le y_j\le z_j\  \text{for}\  j\ne i. $ By the order-preservation, this implies
\beq\label{CMP} \P(X_t^{\bar\mu}\ge X_t^{\mu^i}\ge X_t^{\hat\mu}, t\in [0,T], i=1,2)=1.\end{equation}
Indeed, when $b$ and $\si$ is Lipschitz continuous, by for instance \cite[Theorem 1.1]{HLW} with $\gg=\bar\gg=0$ and $r_0=0$,  \eqref{CMP} follows from \eqref{CMP0} and \eqref{O1}.
Since $b$ and $\si$ are locally Lipschitz continuous and  $X_t^\mu$ is non-explosive for any $\mu\in C_w([0,T];\scr P_\phi)$, we prove \eqref{CMP} by a truncation argument.
Obviously, \eqref{CMP} and $\phi_i'>0$ for $1\le i\le d$ imply
\beq\label{CMP2}    \sum_{i=1}^d \E\big[\phi_i(\{X_t^{\bar\mu}\}_i) - \phi_i(\{X_t^{\hat\mu}\}_i)\big]
 \ge \sum_{i=1}^d \E\big| \phi_i(\{X_t^{\mu^1}\}_i) -  \phi_i(\{X_t^{\mu^2}\}_i)\big|.\end{equation}
Moreover, by \eqref{CMP} and $(A_2)$ we see that $$\xi_t:=  d_\phi(X_t^{\bar\mu}, X_t^{\hat\mu})=\sum_{i=1}^d  \big[\phi_i(\{X_t^{\bar\mu}\}_i) - \phi_i(\{X_t^{\hat\mu}\}_i)\big]$$
satisfies $\xi_t\ge 0$ and
$$\d \xi_t \le (\theta_1+\theta_2-q)  \xi_t\d t +\d M_t$$
for some local martingale $M_t$. As shown in the proof of Lemma \ref{LN1} that for $\ll>2(\theta_1+\theta_2-q)^+$, this implies the contraction of $H$ under the metric $\W_{\phi,\ll}.$

(2)  Next, since $Z_i(x,y)$ is increasing in $(x_j)_{j\ne i}$ and $y$, it is easy to see that conditions  (1) and (2) in \cite[Theorem 1.1]{HLW} holds for $b=\bar b$, and its proof applies also
with $\W_\phi$ replacing $\W_2$ therein, so that the order-preserving property holds.  We omit the details to save space.
Moreover,  since $\scr P_\phi$ is complete under $\W_\phi$, according to the proof of
\cite[Theorem 3.1(2)]{W18}, when $q>\theta_1+\theta_2$ the inequality \eqref{XX1} implies that $P_t^*$ has a unique invariant probability measure $\bar\mu\in \scr P_\phi$ and \eqref{XX2} holds.
Therefore, below we only prove \eqref{XX1}.

(3) To prove    \eqref{XX1}, let $\xi_0,\eta_0$ be $\F_0$-measurable random variable with $\L_{\xi_0}=\mu, \L_{\eta_0}=\nu$ and
\beq\label{X01} \E d_\phi(\xi_0, \eta_0)= \W_\phi(\mu,\nu).\end{equation}
For $x,y\in \R^d$, let $x\lor y=(x_i\lor y_i)_{1\le i\le d}$ and $x\land y=(x_i\land y_i)_{1\le i\le d}.$  Take
\beq\label{X02} X_0=\xi_0\lor\eta_0,\ \ Y_0=\xi_0\land\eta_0.\end{equation}
Let $X_t, Y_t, \xi_t, \eta_t$ solve \eqref{E02} with initial values $X_0, Y_0, \xi_0,\eta_0$ respectively. By the order-preservation, we have
\beq\label{X03}Y_t\le \xi_t\land \eta_t\le \xi_t\lor\eta_t\le X_t,\ \ t\ge 0.\end{equation} Consequently,
$$ d_\phi(X_t,Y_t)=\sum_{i=1}^d \{\phi_i(X_t^i )- \phi_i(Y_t^i)\},\ \ t\ge 0.$$
By It\^o's formula  and applying $(A_1), (A_2)$, we obtain
$$ \d d_\phi(X_t,Y_t)\le \big\{-q d_{\phi}(X_t,Y_t) + \theta_1 d_{\phi}(X_t,Y_t)+ \theta_2 \E d_\phi(X_t,Y_t) \big\}\d t+\d M_t$$
for a local martingale $M_t$.  By a standard argument with Gronwall's lemma,  this implies
$$\E d_\phi (X_t,Y_t)\le \e^{-(q-\theta_1-\theta_2)t} \E d_\phi(X_0,Y_0)=\e^{-(q-\theta_1-\theta_2)t} \E d_\phi(\xi_0,\eta_0),\ \ t\ge 0.$$
Combining this with \eqref{X01} and \eqref{X03}, we arrive at
$$\W_\phi(P_t^*\mu, P_t^*\nu)\le \E d_\phi(\xi_t,\eta_t)\le \E d_\phi(X_t,Y_t)\le \e^{-(q-\theta_1-\theta_2)t} \W_\phi(\mu,\nu),\ \ t\ge 0.$$
Then the proof is finished.

   \paragraph{Acknowledgement.} The author  would like to thank Professor Jian Wang for helpful comments and corrections.

  \end{document}